\documentclass[final,leqno,onefignum,onetabnum]{siamltex1213}
\usepackage{amsmath}
\usepackage{amsfonts}
\usepackage{amssymb}
\usepackage{algorithm,algorithmic}
\usepackage{pst-plot}
\usepackage{graphicx}

\newcommand{\citep}{\cite}
\newcommand{\citet}{\cite}
\def\x{{\mathbf x}}

\def\barg{{ \bar g}}
\def\HH{{\mathbf H}}
\def\z{{\mathbf z}}
\def\1{{\mathbf 1}}

\def\bartheta{{ \bar \theta}}

\def\X{{\mathbf X}}

\def\barh{{ \bar h}}

\def\y{{\mathbf y}}
\def\w{{\mathbf w}}

\def\E{{\mathbb E}}

\def\PPP{{\mathbb P}}
\def\S{{\mathcal S}}

\def\Real{{\mathbb R}}

\def\FF{{\mathcal F}}

\def\argmin{\operatornamewithlimits{arg\,min}}
\def\liminf{\operatornamewithlimits{lim\,inf}}

\def\st{~~\text{s.t.}~~}
\def\defin{\triangleq}

\newcommand{\G}{\mathcal{G}}

\long\def\symbolfootnote[#1]#2{\begingroup\def\thefootnote{\fnsymbol{footnote}}\footnote[#1]{#2}\endgroup} 
\newcommand{\proofstep}[1]{\noindent{\it #1:}~\newline}
 \newcommand{\algorithmicinput}{\textbf{input}}
 \newcommand{\algorithmicoutput}{\textbf{output}}
 \newcommand{\INPUT}{\item[\algorithmicinput]}
 \newcommand{\OUTPUT}{\item[\algorithmicoutput]}
\long\def\proofatend#1\endproofatend{
   \begin{proof}
   #1\end{proof}
}

\title{Incremental Majorization-Minimization Optimization with Application to
   Large-Scale Machine Learning\thanks{This work was partially supported by the
      Gargantua project (program Mastodons - CNRS), the Microsoft
      Research-Inria joint centre, and Agence Nationale
      de la Recherche (MACARON project ANR-14-CE23-0003-01 and the LabEx
      PERSYVAL-Lab ANR-11-LABX-0025). A short version of this work was
      presented at the International Conference of Machine Learning (ICML) in
2013~\cite{mairal17}.}}

\author{Julien Mairal\thanks{Inria, LEAR Team, Laboratoire Jean Kuntzmann, CNRS, Univ. Grenoble Alpes.
655, avenue de l'Europe, 38330 Montbonnot, France. (\email{julien.mairal@inria.fr}).}}

\begin{document}
\maketitle
\slugger{siopt}{xxxx}{xx}{x}{x--x}

\begin{abstract}
   Majorization-minimization algorithms consist of successively minimizing a
   sequence of upper bounds of the objective function. These upper bounds are
   tight at the current estimate, and each iteration monotonically drives the
   objective function downhill. Such a simple principle is widely applicable and
   has been very popular in various scientific fields, especially in signal
   processing and statistics. We propose an incremental majorization-minimization
   scheme for minimizing a large sum of continuous functions, a problem of utmost
   importance in machine learning. We present convergence guarantees for
   non-convex and convex optimization when the upper bounds approximate the
   objective up to a smooth error; we call such upper bounds ``first-order
   surrogate functions''. More precisely, we study asymptotic stationary point
   guarantees for non-convex problems, and for convex ones, we provide convergence
   rates for the expected objective function value. We apply our scheme to
   composite optimization and obtain a new incremental proximal gradient algorithm
   with linear convergence rate for strongly convex functions. Our experiments
   show that our method is competitive with the state of the art for solving
   machine learning problems such as logistic regression when the number of
   training samples is large enough, and we demonstrate its usefulness for sparse
   estimation with non-convex penalties.
\end{abstract}



\pagestyle{myheadings}
\thispagestyle{plain}
\markboth{J. MAIRAL}{INCREMENTAL MAJORIZATION-MINIMIZATION OPTIMIZATION}

\section{Introduction}\label{sec:intro}
The principle of successively minimizing upper bounds of the objective function
is often called \emph{majorization-minimization} \cite{lange2} or
\emph{successive upper-bound minimization}~\cite{razaviyayn}. Each upper bound
is locally tight at the current estimate, and each minimization step decreases
the value of the objective function. Even though this principle does not
provide any theoretical guarantee about the quality of the returned solution,
it has been very popular and widely used because of its simplicity. Various existing
approaches can indeed be interpreted from the majorization-minimization point
of view. This is the case of many gradient-based or proximal
methods~\cite{beck,combette,hale,nesterov,wright}, expectation-maximization (EM) algorithms
in statistics \citep{dempster,neal}, difference-of-convex (DC) programming
\citep{horst}, boosting~\citep{collins,pietra}, some variational Bayes
techniques used in machine learning~\citep{wainwright2}, and the
mean-shift algorithm for finding modes of a
distribution~\cite{tomasi}.  Majorizing surrogates have also been used
successfully in the signal processing literature about sparse estimation
\citep{candes4,daubechies,gasso}, linear inverse problems in image
processing~\cite{ahn,erdogan2}, and matrix factorization \citep{lee2,mairal7}.

In this paper, we are interested in making the
ma\-jo\-ri\-za\-tion-minimization principle scalable for minimizing a large sum
of functions:
\begin{equation}
   \min_{\theta \in \Theta} \left[ f(\theta) \defin \frac{1}{T}\sum_{t=1}^T
   f^t(\theta)\right], \label{eq:prob}
\end{equation}
where the functions $f^t: \Real^p \to \Real$ are continuous, and $\Theta$ is a convex subset of
$\Real^p$. When~$f$ is non-convex, exactly solving~(\ref{eq:prob})
is intractable in general, and when $f$ is also non-smooth, finding a
stationary point of~(\ref{eq:prob}) can be difficult. 
The problem above when $T$ is large can be motivated by machine
learning applications, where~$\theta$ represents some model parameters and
each function $f^t$ measures the adequacy of the parameters~$\theta$ to an
observed data point indexed by~$t$. In this context, minimizing~$f$ amounts to finding
parameters~$\theta$ that explain well some observed data.  In the last few
years, stochastic optimization techniques have become very popular in machine
learning for their empirical ability to deal with a large number $T$ of
training points~\cite{bottou2,duchi2,shalev2,xiao}. Even though these methods
have inherent sublinear convergence rates for convex and strongly convex
problems~\cite{lan,nemirovski}, they typically have a cheap computational cost per
iteration, enabling them to efficiently find an approximate solution.
Recently, incremental algorithms have also been proposed for minimizing finite
sums of functions~\cite{blatt,defazio2,defazio,schmidt2,shalev2}.  At the price of a higher
memory cost than stochastic algorithms, these incremental methods enjoy faster
convergence rates, while also having a cheap per-iteration computational cost.  

Our paper follows this approach: in order to exploit the particular
structure of problem~(\ref{eq:prob}), we propose an incremental scheme whose
cost per iteration is independent of $T$, as soon as the upper bounds of the
objective are appropriately chosen.  We call the resulting scheme ``MISO''
(\emph{Minimization by Incremental Surrogate Optimization}).  We present 
convergence results when the upper bounds are chosen among the class of
``first-order surrogate functions'', which approximate the objective function
up to a smooth error---that is, differentiable with a Lipschitz continuous
gradient.  For non-convex problems, we obtain almost sure convergence and
asymptotic stationary point guarantees. In addition, when assuming the surrogates
to be strongly convex, we provide convergence rates for the expected value of
the objective function. Remarkably, the convergence rate of
MISO is linear for minimizing strongly convex composite objective functions, a
property shared with two other incremental algorithms for smooth and composite
convex optimization: the \emph{stochastic average gradient} method (SAG) of Schmidt,
Le Roux and Bach~\citet{schmidt2}, and the \emph{stochastic dual coordinate
ascent} method (SDCA) of Shalev-Schwartz and Zhang~\citet{shalev2}.  Our scheme
MISO is inspired in part by these two works, but yields different update rules than
SAG or SDCA, and is also appropriate for non-convex optimization problems.

In the experimental section of this paper, we show that MISO can be useful for
solving large-scale machine learning problems, and that it matches cutting-edge solvers for large-scale logistic regression
\citep{beck,schmidt2}. Then, we show that our approach
provides an effective incremental DC programming algorithm, which we apply to
sparse estimation problems with nonconvex penalties~\cite{candes4}. 

The paper is organized as follows: Section~\ref{sec:generic} introduces the
majorization-mi\-ni\-mi\-za\-tion principle with first-order surrogate
functions. Section~\ref{sec:incremental} is devoted to our incremental scheme
MISO. Section~\ref{sec:exp} presents some numerical experiments, and
Section~\ref{sec:ccl} concludes the paper. Some basic definitions are given in
Appendix~\ref{appendix:background}.

\section{Majorization-minimization with first-order surrogate functions}\label{sec:generic}
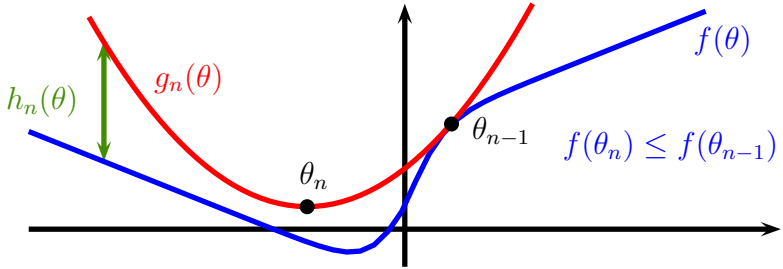
\begin{figure}[hbtp]
   \centering
   \psset{yunit=1,xunit=1}
   \begin{pspicture}(-5,-0.5)(5,4)
      \psline[linewidth=2pt]{->}(-5,0)(5,0)
      \psline[linewidth=2pt]{->}(0,-0.5)(0,3)
      \psline[linecolor=darkgreen, linewidth=2pt]{<->}(-4,0.9)(-4,2.5)
      \psplot[linecolor=blue, linewidth=2pt]{-5}{4}{
         2.0 1.0 2.7183 x -4.0 mul exp 1.0 add div x abs 0.2 mul add mul -0.7 add
      }
      \psplot[linecolor=red, linewidth=2pt]{-4.2}{1.7}{
         x 1.3 add x 1.3 add mul 0.3 mul 0.3 add
      } 
      \put(3.84,2.4){\color{blue}  \large ${ f(\theta)}$}
      \put(-3.33,1.9){\color{red}  \large ${ g_n(\theta)}$}
      \put(-5.3,1.6){\color{darkgreen}  \large ${ h_n(\theta)}$}
      \psdots[dotsize=0.2](0.62,1.40)(-1.3,0.3)
      \put(0.9,1.2) {\large $\theta_{n-1}$}
      \put(-1.4,0.6) {\large $\theta_{n}$}
      \put(2.1,1.0)  {\color{blue} \large $f(\theta_n) \leq f(\theta_{n-1})$}
   \end{pspicture}
   \caption{Illustration of the basic majorization-minimization principle. We compute a
      surrogate~$g_n$ of~$f$ near the current
      estimate~$\theta_{n-1}$. The new estimate~$\theta_n$ is a minimizer
      of~$g_n$. The function~$h_n=g_n-f$ is the approximation error that is made when
   replacing $f$ by~$g_n$.}\label{fig:mm} 
\end{figure}

In this section, we present the generic majorization-minimization scheme for
minimizing a function~$f$ without exploiting its structure---that is,
without using the fact that $f$ is a sum of functions.
We describe the procedure in Algorithm~\ref{alg:generic_batch} and illustrate its principle 
in Figure~\ref{fig:mm}. At iteration~$n$, the estimate~$\theta_n$ is obtained by
minimizing a surrogate function~$g_n$ of~$f$. When $g_n$
uniformly majorizes~$f$ and when $g_n(\theta_{n-1})=f(\theta_{n-1})$, it is
clear that the objective function value monotonically decreases. 

\begin{algorithm}[hbtp]
   \caption{Basic majorization-minimization scheme.}\label{alg:generic_batch}
   \begin{algorithmic}[1]
      \INPUT $\theta_0 \in \Theta$ (initial estimate); $N$ (number of iterations).
      \FOR{ $n=1,\ldots,N$}
      \STATE Compute a surrogate function $g_n$ of $f$ near $\theta_{n-1}$;
      \STATE Minimize the surrogate and update the solution:
      $\theta_n \in \argmin_{\theta \in \Theta} g_n(\theta).$
   \ENDFOR
   \OUTPUT $\theta_{N}$ (final estimate);
\end{algorithmic}
\end{algorithm}

\definecolorset{rgb}{}{}{darkgreen,0.26,0.55,0}

For this approach to be effective, intuition tells us that we need functions~$g_n$ that
are easy to minimize and that approximate well the objective~$f$. Therefore,
we measure the quality of the approximation through the smoothness of the
error $h_n \defin g_n-f$, which is a key quantity arising in the
convergence analysis. Specifically, we require~$h_n$ to be $L$-smooth for
some constant~$L>0$ in the following sense:

\begin{definition}[{\rm $L$-smooth functions}]\label{def:lsmooth}
   A function $f: \Real^p \to \Real$ is called $L$-smooth when it is differentiable and when its gradient $\nabla f$ is $L$-Lipschitz continuous.
\end{definition}

With this definition in hand, we now introduce the class of ``first-order
surrogate functions'', which will be shown to have good enough properties for
analyzing the convergence of Algorithm~\ref{alg:generic_batch} and other
variants.

\begin{definition}[{\rm First-order surrogate functions}]\label{def:surrogate_batch}
   A function $g: \Real^p \to \Real$ is a first-order surrogate function of $f$
   near~$\kappa$ in~$\Theta$ when
   \begin{romannum}
      \item $g(\theta') \geq f(\theta')$ for all
         minimizers $\theta'$ of $g$ over $\Theta$.
         When the more general condition $g \geq f$ holds, we say that $g$ is a \emph{majorizing} surrogate;
      \item the approximation error $h\defin g-f$ is $L$-smooth, $h(\kappa)=0$,
         and $\nabla h(\kappa)=0$. 
   \end{romannum}
   We denote by~$\S_{L}(f,\kappa)$ the set of first-order surrogate functions and by~$\S_{L,\rho}(f,\kappa) \subset \S_{L}(f,\kappa)$ the subset of $\rho$-strongly convex surrogates.
\end{definition}

First-order surrogates are interesting because their approximation error---the
difference between the surrogate and the objective---can be easily controlled.
This is formally stated in the next lemma, which is a building block of our
analysis:
\begin{lemma}[{\rm Basic properties of first-order surrogate functions}]\label{lemma:basic}
   Let~$g$ be a surrogate function in $\S_{L}(f,\kappa)$ for some $\kappa$ in $\Theta$. Define the approximation error $h\defin g-f$,
   and let $\theta'$ be a minimizer of $g$ over $\Theta$. Then, for all~$\theta$ in~$\Theta$,
   \begin{itemize}
      \item $|h(\theta)| \leq \frac{L}{2}\|\theta-\kappa\|_2^2$;
      \item $f(\theta') \leq f(\theta) + \frac{L}{2}\|\theta-\kappa\|_2^2$.
   \end{itemize}
   Assume that $g$ is $\rho$-strongly convex, i,e., $g$ is in $\S_{L,\rho}(f,\kappa)$. Then, for all $\theta$ in $\Theta$, 
   \begin{itemize}
      \item $f(\theta') + \frac{\rho}{2}\|\theta'-\theta\|_2^2 \leq f(\theta) + \frac{L}{2}\|\theta-\kappa\|_2^2$.
   \end{itemize}
\end{lemma}

\proofatend
The first inequality is a direct application of a classical result (Lemma 1.2.3
of~\cite{nesterov4}) on quadratic upper bounds for $L$-smooth functions, when
noticing that $h(\kappa)=0$ and $\nabla h(\kappa)=0$. Then, for all $\theta$
in $\Theta$, we have $f(\theta') \leq g(\theta') \leq g(\theta) =
f(\theta)+h(\theta),$ and we obtain the second inequality from the first one.

When $g$ is $\rho$-strongly convex, we use the following classical lower bound~(see \cite{nesterov}):
$$ g(\theta') +
\frac{\rho}{2}\|\theta-\theta'\|_2^2 \leq g(\theta).$$
Since $f(\theta') \leq g(\theta')$ by Definition~\ref{def:surrogate_batch} and $g(\theta)=f(\theta)+h(\theta)$,
the third inequality follows from the first one.
\endproofatend

We now proceed with a convergence analysis including four main results
regarding Algorithm~\ref{alg:generic_batch} with first-order surrogate
functions~$g_n$. More precisely, we show in Section~\ref{subsec:nonconvex}
that, under simple assumptions, the sequence of iterates asymptotically
satisfies a stationary point condition. Then, we present a similar result with
relaxed assumptions on the surrogates~$g_n$ when $f$ is a composition of two
functions, which occur in practical situations as shown in
Section~\ref{subsec:surrogates}. Finally, we present non-asymptotic convergence
rates when~$f$ is convex in Section~\ref{subsec:convex}. By adapting 
convergence proofs of proximal gradient methods~\citep{nesterov} to our more
general setting, we recover classical sublinear rates $O(1/n)$
and linear convergence rates for strongly convex problems.

\subsection{Non-convex convergence analysis}\label{subsec:nonconvex}
For general non-convex problems, proving convergence to a global (or local)
minimum is impossible in general, and classical analysis studies instead asymptotic
stationary point conditions~(see, {\it e.g.}, \cite{bertsekas}).
To do so, we make the following mild assumption when $f$ is non-convex:
\def\assumpBasic{{\rm (A)}}
\begin{itemize}
   \item[\assumpBasic] $f$ is bounded below and for all $\theta,\theta'$ in $\Theta$, the
      directional derivative $\nabla f(\theta,\theta'-\theta)$ of $f$ at~$\theta$ in the direction
      $\theta'-\theta$ exists.
\end{itemize}
The definitions of directional derivatives and stationary points are provided in
Appendix~\ref{appendix:background}. A necessary first-order
condition for~$\theta$ to be a local minimum of $f$ is to have $\nabla
f(\theta,\theta'\!-\!\theta) \geq 0$ for all~$\theta'$ in~$\Theta$ (see,
{\it e.g.}, \cite{borwein}).  In other words, there is no feasible descent
direction~$\theta'\!-\!\theta$ and $\theta$ is a stationary point. Thus, we consider
the following condition for assessing the quality of a
sequence $(\theta_n)_{n \geq 0}$ for non-convex problems:
\begin{definition}[{\rm Asymptotic stationary point}]
   Under assumption~\assumpBasic, a sequence $(\theta_n)_{n \geq 0}$ satisfies the asymptotic stationary point condition if
   \begin{equation}
      \liminf_{n \to +\infty} \inf_{\theta \in \Theta} \frac{\nabla
      f(\theta_{n},\theta-\theta_n)}{ \|\theta-\theta_n\|_2} \geq 0. \label{eq:stationary}
   \end{equation}
   Note that if $f$ is differentiable on $\Real^p$ and $\Theta = \Real^p$,
   $\nabla f(\theta_{n},\theta-\theta_n) =\nabla
   f(\theta_{n})^\top(\theta-\theta_n)$,
   and the condition~(\ref{eq:stationary}) implies that the sequence $(\nabla f(\theta_n))_{n \geq 0}$ converges to $0$.
\end{definition}
As noted, we recover the classical definition of critical points for the smooth
unconstrained case. We now give a first convergence result about
Algorithm~\ref{alg:generic_batch}.

\begin{proposition}[\rm Non-convex analysis for Algorithm~\ref{alg:generic_batch}]\label{prop:conv1}
   Assume that~\assumpBasic~holds and that the surrogates $g_n$ from Algorithm~\ref{alg:generic_batch} are in $\S_L(f,\theta_{n-1})$ and
   are either majorizing $f$ or strongly convex.  Then, $\!(f(\theta_n))_{n \geq 0}$
   monotonically decreases, and~$(\theta_n)_{n\geq 0}$ satisfies the asymptotic
   stationary point condition.
\end{proposition}
\proofatend
The fact that $(f(\theta_n))_{n \geq 0}$ is non-increasing and convergent because bounded below is clear: for all $n\geq 1$,
$f(\theta_n) \leq g_n(\theta_n) \leq g_n(\theta_{n-1}) = f(\theta_{n-1})$,
where the first inequality and the last equality are obtained from
Definition~\ref{def:surrogate_batch}. The second inequality comes from the
definition of $\theta_n$.

Let us now denote by $f^\star$ the limit of the sequence $(f(\theta_n))_{n \geq
1}$ and by $h_n\defin g_n-f$ the approximation error function at iteration~$n$,
which is $L$-smooth by Definition~\ref{def:surrogate_batch} and such
that~$h_n(\theta_n) \geq 0$.  Then,
$h_n(\theta_n) = g_n(\theta_n) - f(\theta_n) \leq f(\theta_{n-1}) - f(\theta_n)$, and
\begin{displaymath}
   \sum_{n=1}^{\infty} h_n(\theta_n) \leq f(\theta_0) - f^\star.
\end{displaymath}
Thus, the non-negative sequence $(h_n(\theta_n))_{n \geq 0}$ necessarily converges to zero.
Then, we have two possibilities (according to the assumptions made in the proposition).
\begin{itemize}

   \item If the functions $g_n$ are majorizing~$f$, we define $\theta' = \theta_{n}-\frac{1}{L}\nabla h_n(\theta_{n})$, and 
      we use the following classical inequality for $L$-smooth functions~\cite{nesterov4}:
      \begin{displaymath}
         h_n(\theta') \leq h_n(\theta_n) - \frac{1}{2L}\|\nabla h_n(\theta_n)\|_2^2.
      \end{displaymath}
      Therefore, we may use the fact that $h_n(\theta') \geq 0$ because $g_n \geq f$, and 
      \begin{displaymath}
         \|\nabla h_n(\theta_n)\|_2^2 \leq 2L (h_n(\theta_n)-h_n(\theta')) \leq 2L h_n(\theta_n) \underset{n \to +\infty}{\longrightarrow} 0.
      \end{displaymath}

   \item If instead the functions $g_n$ are $\rho$-strongly convex, the last inequality of Lemma~\ref{lemma:basic} with $\kappa = \theta = \theta_{n-1}$ and~$\theta' = \theta_n$ gives us
      \begin{displaymath} 
         \frac{\rho}{2}\|\theta_n-\theta_{n-1}\|_2^2  \leq f(\theta_{n-1})-f(\theta_n).
      \end{displaymath}
      By summing over $n$, we obtain that $\|\theta_n-\theta_{n-1}\|_2^2$ converges to zero, and
      \begin{displaymath} 
         \|\nabla h_n(\theta_n)\|_2 =    \|\nabla h_n(\theta_n)-\nabla h_n(\theta_{n-1})\|_2 \leq L\|\theta_n-\theta_{n-1}\|_2 \underset{n \to +\infty}{\longrightarrow} 0,
      \end{displaymath}
      since $\nabla h_n(\theta_{n-1})=0$ according to Definition~\ref{def:surrogate_batch}.
\end{itemize}

We now consider the directional derivative of $f$ at $\theta_n$ and a direction $\theta-\theta_n$, where $n \geq 1$ and $\theta$ is in~$\Theta$,
\begin{displaymath}
   \nabla f(\theta_{n},\theta-\theta_n) = \nabla g_n(\theta_{n},\theta-\theta_n) - \nabla h_n(\theta_n)^\top (\theta-\theta_n).
\end{displaymath}
Note that $\theta_n$ minimizes $g_n$ on $\Theta$ and therefore $\nabla g_n(\theta_{n},\theta-\theta_n) \geq 0$. Therefore, 
\begin{displaymath}
   \nabla f(\theta_{n},\theta-\theta_n)  \geq - \|\nabla h_n(\theta_n)\|_2 \|\theta-\theta_n\|_2,
\end{displaymath}
by Cauchy-Schwarz's inequality. By minimizing over $\theta$ and taking the infimum limit, we finally obtain
\begin{displaymath}
   \liminf_{n \to +\infty} \inf_{\theta \in \Theta} \frac{\nabla f(\theta_{n},\theta-\theta_n)}{ \|\theta-\theta_n\|_2} \geq - \lim_{n\to+\infty} \|\nabla h_n(\theta_n)\|_2 = 0.
\end{displaymath}
\endproofatend

This proposition provides convergence guarantees for a large class of existing
algorithms, including cases where $f$ is non-smooth. In the next proposition,
we relax some of the assumptions for objective functions that are compositions $f=f' \circ e$,
where $\circ$ is the composition operator. In other words, $f(\theta) =
f'(e(\theta))$ for all $\theta$ in~$\Real^p$.
\begin{proposition}[\rm Non-convex analysis for Algorithm~\ref{alg:generic_batch} - composition]\label{prop:conv1_separable}
   Assume that~\assumpBasic~holds and that the function $f$ is a composition
   $f = f' \circ e$, where $e: \Real^p \to \Real^d$ is $C$-Lipschitz continuous
   for some constant $C>0$, and $f': \Real^d \to \Real$. Assume that the function~$g_n$
   in Algorithm~\ref{alg:generic_batch} is defined as $g_{n}\defin
   g_{n}' \circ e$, where $g_{n}'$ is a majorizing surrogate in
   $\S_L(f',e(\theta_{n-1}))$.  Then, the conclusions of
   Proposition~\ref{prop:conv1} hold.
\end{proposition}
\proofatend
We follow the same steps as the proof of Proposition~\ref{prop:conv1}.
First, it is easy to show that $(f(\theta_n))_{n \geq 0}$ monotonically decreases and
that $h_n(\theta_n) \defin g_n(\theta_n)-f(\theta_n)$ converges to zero when
$n$ grows to infinity.  Note that since we have made the assumptions that $g_n = g_n' \circ e$ and that $f = f' \circ e$,  
the function $h_n\defin  g_n - f$ can be written as $h_n = h_n' \circ e$,
where $h_n'\defin g_n'- f'$ is $L$-smooth.  Proceeding as in the proof of
Proposition~\ref{prop:conv1}, we can  show that $\|\nabla
h_n'(e(\theta_n))\|_2$ converges to zero.

Let us now fix $n \geq 1$ and consider $\delta$ such that $\theta_n+\delta$ is in~$\Theta$. We have
\begin{displaymath}
   h_n(\theta_n+\delta) = h_n'( e(\theta_n+\delta)) = h_n'( e(\theta_n) + \|\delta\|_2 \z), 
\end{displaymath}
where $\z$ is a vector whose $\ell_2$-norm is bounded by a universal
constant~$C > 0$ because the function~$e$ is Lipschitz continuous.
Since $h_n'$ is $L$-smooth, we also have
\begin{displaymath}
   h_n(\theta_n+\delta) = h_n'( e(\theta_n) + \|\delta\|_2 \z) = h_n( \theta_n) + \|\delta\|_2 \nabla h_n^{\prime}(e(\theta_n))^\top \z + O(\|\delta\|_2^2).
\end{displaymath}
Plugging this simple relation with $\delta=t(\theta-\theta_n)$, for some $0 < t < 1$ and $\theta$ in~$\Theta$, into the definition of the
directional derivative $\nabla h_n(\theta_n, \theta-\theta_n)$, we obtain the
relation
\begin{displaymath}
   |\nabla h_n(\theta_n, \theta-\theta_n)| \leq C\| \nabla h_n^{\prime}(e(\theta_n))\|_2{\|\theta-\theta_n\|_2},
\end{displaymath}
and since $\nabla f(\theta_{n},\theta-\theta_n) = \nabla g_n(\theta_{n},\theta-\theta_n) - \nabla h_n(\theta_n,\theta-\theta_n)$, and $\nabla g_n(\theta_{n},\theta-\theta_n) \geq 0$, 
\begin{displaymath}
   \liminf_{n \to +\infty} \inf_{\theta \in \Theta} \frac{\nabla f(\theta_{n},\theta-\theta_n)}{ \|\theta-\theta_n\|_2} \geq -C \lim_{n \to +\infty} \| \nabla h_n^{\prime}(e(\theta_n))\|_2 = 0.
\end{displaymath}
\endproofatend

In this proposition, $g_n$ is an upper bound of~$f=f' \circ e$, where the
part~$e$ is Lipschitz continuous but $g_n-f$ is not~$L$-smooth. This extension
of Proposition~\ref{prop:conv1} is useful since it provides convergence results
for classical approaches that will be described later in
Section~\ref{subsec:surrogates}.  Note that convergence results for non-convex
problems are by nature weak, and our non-convex analysis does not provide any
convergence rate.  This is not the case when~$f$ is convex, as shown in the
next section.

\subsection{Convex analysis}\label{subsec:convex}
The next proposition is based on a proof technique from
Nesterov~\citet{nesterov}, which was originally designed for the proximal gradient
method. By adapting it, we obtain the same
convergence rates as in~\citep{nesterov}.

\begin{proposition}[\rm Convex analysis for $\S_L(f,\kappa)$]\label{prop:conv2}
   Assume that $f$ is convex, bounded below, and that there exists a constant $R >0$ such that
   \begin{equation}
      \|\theta-\theta^\star\|_2 \leq R ~~~\text{for all}~ \theta \in \Theta \st f(\theta) \leq f(\theta_0),\label{eq:bounded}
   \end{equation}
   where $\theta^\star$ is a minimizer of $f$ on $\Theta$. 
   When the functions~$g_n$ in Algorithm~\ref{alg:generic_batch} are in $\S_{L}(f,\theta_{n-1})$, we have for all $n \geq 1$, 
   $$f(\theta_n)-f^\star \leq \frac{2LR^2}{n+2},$$
   where $f^\star \defin f(\theta^\star)$. 
   Assume now that $f$ is $\mu$-strongly convex. Regardless of condition~(\ref{eq:bounded}), we have for all $n \geq 1$,
   \begin{displaymath}
      f(\theta_n) - f^\star  \leq \beta^n( f(\theta_0)-f^\star),
   \end{displaymath}
   where $\beta \defin \frac{L}{\mu}$ if $\mu >2L$ or $\beta \defin \left(1-\frac{\mu}{4L}\right)$ otherwise.
\end{proposition}
\proofatend 
We successively prove the two parts of the proposition.\\
\proofstep{Non-strongly convex case}
Let us consider the function $h_n \defin g_n - f$ at iteration $n \geq 1$. By Lemma~\ref{lemma:basic}, 
\begin{displaymath}
   f(\theta_n) \leq \min_{\theta \in \Theta} \left[ f(\theta) + \frac{L}{2} \|\theta-\theta_{n-1}\|_2^2\right].
\end{displaymath}
Then, following a similar proof technique as Nesterov in~\cite{nesterov}, 
\begin{equation}
   \begin{split}
   f(\theta_n) & \leq \min_{\alpha \in [0,1]} \left[ f(\alpha\theta^\star+(1-\alpha)\theta_{n-1}) + \frac{L\alpha^2}{2} \|\theta^\star-\theta_{n-1}\|_2^2\right] \\
                  & \leq \min_{\alpha \in [0,1]}  \left[ \alpha f(\theta^\star) +(1-\alpha)f(\theta_{n-1}) + \frac{L\alpha^2}{2} \|\theta^\star-\theta_{n-1}\|_2^2\right],
   \end{split} \label{eq:tmp_rate2}
\end{equation}
where the minimization over $\Theta$ is replaced by a
minimization over the line segment $\alpha\theta^\star+(1-\alpha)\theta_{n-1} :
\alpha \in [0,1]$. Since the sequence $(f(\theta_n))_{n \geq 0}$ is
monotonically decreasing we may use the bounded level set assumption and we obtain
\begin{equation}
   f(\theta_n) - f^\star \leq \min_{\alpha \in [0,1]} \left[ (1-\alpha)(f(\theta_{n-1}) - f^\star) + \frac{LR^2\alpha^2}{2}\right].\label{eq:tmp_rate3}
\end{equation}
To simplify, we introduce the notation~$r_n \defin f(\theta_n)- f^\star$, and we consider two cases:
\begin{itemize}
   \item {\it first case:} if $r_{n-1} \geq LR^2$, then the optimal value $\alpha^\star$ in~(\ref{eq:tmp_rate3}) is~$1$ and we consequently have $r_n \leq \frac{LR^2}{2}$;
   \item {\it second case:} otherwise $\alpha^\star =\frac{r_{n-1}}{LR^2}$ and $r_n \leq r_{n-1}\left(1-\frac{r_{n-1}}{2LR^2}\right)$.
      Thus, $r_n^{-1} \geq r_{n-1}^{-1} \left(1-\frac{r_{n-1}}{2LR^2}\right)^{-1} \geq r_{n-1}^{-1} + \frac{1}{2LR^2}$, where the second inequality comes from the convexity inequality $(1-x)^{-1} \geq 1+x$ for $x \in (0,1)$.
\end{itemize}
We now apply recursively the previous inequalities, starting with $n=1$. If $r_0 \geq LR^2$, we are in the first case and then $r_1 \leq \frac{LR^2}{2}$; Then, we will subsequently be in the second
case for all $n \geq 2$ and thus $r_n^{-1} \geq  r_1^{-1} + \frac{n-1}{2LR^2} \geq \frac{n+3}{2LR^2}$. Otherwise, if $r_0 < LR^2$, we are always in the second case and $r_n^{-1} \geq r_0^{-1} + \frac{n}{2LR^2} \geq \frac{n+2}{2LR^2}$, which is sufficient to obtain the first part of the proposition.

\proofstep{$\mu$-strongly convex case}
Let us now assume that $f$ is $\mu$-strongly convex, and let us drop the bounded level
sets assumption.  The proof again follows~\citet{nesterov} for computing the
convergence rate of proximal gradient methods. We start from~(\ref{eq:tmp_rate2}).
We use the strong convexity of~$f$ which implies that $f(\theta_{n-1}) \geq f^\star + \frac{\mu}{2}\|\theta_{n-1}-\theta^\star\|_2^2$, and we obtain
\begin{displaymath}
   f(\theta_n) - f^\star \leq \left(\min_{\alpha \in [0,1]}  1-\alpha + \frac{L\alpha^2}{\mu}\right)(f(\theta_{n-1})-f^\star).
\end{displaymath}
At this point, it is easy to show that 
if $\mu \geq 2L$, the previous binomial is minimized for $\alpha^\star = 1$, and
if $\mu \leq 2L$, then we have $\alpha^\star = \frac{\mu}{2L}$.
This yields the desired result.

\endproofatend

The result of Proposition~\ref{prop:conv2} is interesting because it does not
make any strong assumption about the surrogate functions, except the ones from
Definition~\ref{def:surrogate_batch}.  The next proposition shows that slightly
better rates can be obtained with additional strong convexity assumptions.
\begin{proposition}[\rm Convex analysis for $\S_{L,\rho}(f,\kappa)$]\label{prop:conva}
   Assume that $f$ is convex, bounded below, and let $\theta^\star$ be a minimizer of~$f$ on~$\Theta$.
   When the surrogates~$g_n$ of Algorithm~\ref{alg:generic_batch} are in $\S_{L,\rho}(f,\theta_{n-1})$ with $\rho \geq L$, we have for all $n \geq 1$,
   \begin{displaymath}
      f(\theta_n) - f^\star  \leq  \frac{L\|\theta_{0}-\theta^\star\|_2^2}{2n},\label{eq:rate1}
   \end{displaymath}
   where $f^\star \defin f(\theta^\star)$.
   When $f$ is $\mu$-strongly convex, we have for all $n\geq 1$,
   \begin{displaymath}
                f(\theta_n) - f^\star  \leq \left(\frac{L}{\rho+\mu}\right)^{n-1}\frac{L\|\theta_0-\theta^\star\|_2^2}{2}.
       \end{displaymath} 
   \end{proposition}
   \proofatend 
   As before, we successively prove the two parts of the proposition.

   \proofstep{Non-strongly convex case}
   From Lemma~\ref{lemma:basic} (with $g=g_n$, $\kappa=\theta_{n-1}$, $\theta'=\theta_n$, $\theta=\theta^\star$), we have for all $n\geq 1$,
   \begin{equation}
      f(\theta_n) - f(\theta^\star) \leq \frac{L}{2}\|\theta_{n-1}-\theta^\star\|_2^2 - \frac{\rho}{2}\|\theta_n-\theta^\star\|_2^2\leq \frac{L}{2}\|\theta_{n-1}-\theta^\star\|_2^2 - \frac{L}{2}\|\theta_n-\theta^\star\|_2^2.\label{eq:rate1_proof}
   \end{equation}
   After summation, 
   $$ n(f(\theta_n) - f(\theta^\star)) \leq \sum_{k=1}^n (f(\theta_k) - f(\theta^\star))\leq \frac{L}{2}(\|\theta_{0}-\theta^\star\|_2^2-\|\theta_{n}-\theta^\star\|_2^2) \leq \frac{L\|\theta_{0}-\theta^\star\|_2^2}{2},$$
   where the first inequality comes from the inequalities $f(\theta_k) \geq f(\theta_n)$ for all $k \leq n$. This is sufficient to prove the first part. Note that proving convergence rates for first-order methods by finding telescopic sums is a classical technique~(see, {\it e.g.},\cite{beck}).

   \proofstep{$\mu$-strongly convex case}
   Let us now assume that $f$ is $\mu$-strongly convex. The strong convexity implies that $f(\theta_n)-f^\star \geq \frac{\mu}{2}\|\theta_n-\theta^\star\|_2^2$ for all $n$. Combined with~(\ref{eq:rate1_proof}), this yields
   $$
   \frac{\mu+\rho}{2}\|\theta_n-\theta^\star\|_2^2 \leq \frac{L}{2}\|\theta_{n-1}-\theta^\star\|_2^2,
   $$
   and thus
   $$
   f(\theta_n) - f(\theta^\star) \leq \frac{L}{2}\|\theta_{n-1}-\theta^\star\|_2^2 \leq \left(\frac{L}{\rho+\mu}\right)^{n-1}\frac{L\|\theta_0-\theta^\star\|_2^2}{2}.
   $$
\endproofatend

Even though the constants obtained in the rates of Proposition~\ref{prop:conva} are
slightly better than the ones of Proposition~\ref{prop:conv2}, 
the condition $g_n$ in $\S_{L,\rho}(f,\kappa)$ with $\rho \geq L$ is much stronger 
than the simple assumption that~$g_n$ is in $\S_L(f,\kappa)$.
It can indeed be shown that~$f$ is necessarily 
$(\rho\!-\!L)$-strongly convex if $\rho \!>\! L$, and convex if $\rho\!=\!L$.
In the next section, we give some examples where such a condition holds.
\subsection{Examples of first-order surrogate functions} \label{subsec:surrogates}
We now present practical first-order surrogate functions and links between
Algorithm~\ref{alg:generic_batch} and existing approaches. Even though our
generic analysis does not always bring new results for each specific case, its
main asset is to provide a unique theoretical treatment to all of
them.

\subsubsection{Lipschitz gradient surrogates}\label{subsubsec:gradient}
When $f$ is~$L$-smooth, it is natural to consider the following surrogate:
\begin{displaymath}
    g: \theta \mapsto f(\kappa) + \nabla f(\kappa)^\top (\theta-\kappa) + \frac{L}{2}\|\theta-\kappa\|_2^2.
\end{displaymath}
The function~$g$ is an upper bound of $f$, which is a classical
result~\cite{nesterov4}. It is then easy to see that~$g$ is~$L$-strongly convex
and $L$-smooth.  As a consequence, the difference~$g-f$ is~$2L$-smooth (as a
sum of two~$L$-smooth functions), and thus $g$ is in $\S_{2L,L}(f,\kappa)$.

When $f$ is convex, it is also possible to show by using
Lemma~\ref{lemma:convexerror} that $g$ is in fact in $\S_{L,L}(f,\kappa)$, and
when $f$ is $\mu$-strongly convex, $g$ is in $\S_{L-\mu,L}(f,\kappa)$. We
remark that mi\-ni\-mi\-zing~$g$ amounts to performing a gradient
descent step: $\theta' \leftarrow \kappa - \frac{1}{L}\nabla f(\kappa)$.

\subsubsection{Proximal gradient surrogates}\label{subsubsec:proximal}
Let us now consider a composite optimization problem, meaning that $f$ splits
into two parts $f = f_1 + f_2$, where $f_1$ is $L$-smooth. Then, a natural surrogate of~$f$ 
is the following function:
\begin{displaymath}
    g: \theta \mapsto f_1(\kappa) + \nabla f_1(\kappa)^\top (\theta-\kappa) + \frac{L}{2}\|\theta-\kappa\|_2^2 + f_2(\theta).
\end{displaymath}
The function~$g$ majorizes~$f$ and the approximation error $g-f$ is the same as in Section~\ref{subsubsec:gradient}. Thus,
$g$ in in~$\S_{2L}(f,\kappa)$ or in $\S_{2L,L}(f,\kappa)$ when $f_2$ is convex. Moreover,
\begin{itemize}
   \item when $f_1$ is convex, $g$ is in $\S_{L}(f,\kappa)$. If $f_2$ is also convex, $g$ is in $\S_{L,L}(f,\kappa)$;
   \item when $f_1$ is $\mu$-strongly convex, $g$ is in $\S_{L-\mu}(f,\kappa)$. If $f_2$ is also convex, $g$ is in $\S_{L-\mu,L}(f,\kappa)$.
\end{itemize}
Minimizing~$g$ amounts to performing one step of the proximal gradient algorithm~\cite{beck,nesterov,wright}.
It is indeed easy to show that the minimum $\theta'$ of $g$---assuming it is
unique---can be equivalently obtained as follows:
 \begin{displaymath}
    \theta'  = \argmin_{\theta \in \Theta} \left[\frac{1}{2}\left\| \theta - \left(\kappa -\frac{1}{L}\nabla f_1(\kappa)\right)\right\|_2^2 + \frac{1}{L}f_2(\theta)\right],
 \end{displaymath}
which is often written under the form $\theta' =
\text{Prox}_{f_2/L}[\kappa-(1/L)\nabla f_1(\kappa)]$, where ``$\text{Prox}$'' is
called the ``proximal operator''~\cite{moreau}. In some cases, the proximal operator can be
computed efficiently in closed form, for example when $f_2$ is the~$\ell_1$-norm; it yields
the iterative soft-thresholding algorithm for sparse
estimation~\citep{daubechies}.  For a review of proximal operators and their
computations, we refer the reader to~\citep{bach8,combettes2005signal}.
 
\subsubsection{Linearizing concave functions and DC programming} \label{subsec:dc}
Assume that $f = f_1 + f_2$, where $f_2$ is concave and~$L$-smooth. Then, the following function $g$ is a majorizing
surrogate in $\S_{L}(f,\kappa)$:
\begin{displaymath}
   g: \theta \mapsto f_1(\theta) + f_2(\kappa) + \nabla f_2(\kappa)^\top (\theta-\kappa).
\end{displaymath}
Such a surrogate appears in DC (difference of convex) programming~\cite{horst}. When $f_1$ is convex,
$f$ is indeed the difference of two convex functions.
It is also used in sparse estimation for dealing with some non-convex
penalties~\cite{bach8}. For example, consider a cost function of the form
$\theta \mapsto f_1(\theta) + \lambda \sum_{j=1}^p \log(|\theta[j]| +
\varepsilon)$, where $\theta[j]$ is the $j$-th entry in~$\theta$. Even though the functions $\theta \mapsto \log(|\theta[j]|
+\varepsilon)$ are not differentiable, they can be written as the
composition of a concave smooth function $u \mapsto \log(u + \varepsilon)$ on~$\Real^+$, and a Lipschitz function $\theta \mapsto |\theta[j]|$. By upper-bounding the logarithm function
by its linear approximation, Proposition~\ref{prop:conv1_separable} justifies the following surrogate:
\begin{equation}
   g: \theta \mapsto f_1(\theta) + \lambda \sum_{j=1}^p \log(|\kappa[j]| + \varepsilon) + \lambda\sum_{j=1}^p \frac{|\theta[j]|-|\kappa[j]|}{|\kappa[j]|+\varepsilon}, \label{eq:upperbounddc}
\end{equation}
and minimizing $g$ amounts to performing one step of a reweighted-$\ell_1$
algorithm (see~\citet{candes4} and references therein).  Similarly, other penalty
functions are adapted to this framework. For instance, the logarithm can be replaced by any smooth concave non-decreasing function, or 
group-sparsity penalties~\cite{turlach,yuan} can be used, such as $\theta \mapsto \sum_{g \in \G}
\log(\|\theta_g\|_2 + \varepsilon)$, where $\G$ is a partition of
$\{1,\ldots,p\}$ and $\theta_g$ records the entries of~$\theta$ corresponding
to the set $g$.  Proposition~\ref{prop:conv1_separable} indeed applies to this setting.

\subsubsection{Variational surrogates}\label{subsubsec:variational}
Let us now consider a
real-valued function $f$ defined on $\Real^{p_1} \times \Real^{p_2}$. Let
$\Theta_1 \subseteq \Real^{p_1}$ and $\Theta_2 \subseteq \Real^{p_2}$ be two
convex sets. Minimizing~$f$ over~$\Theta_1 \times \Theta_2$ is equivalent to
minimizing the function~$\tilde{f}$ over~$\Theta_1$ defined as $\tilde{f}(\theta_1) \defin \mapsto \min_{\theta_2
\in \Theta_2} f(\theta_1,\theta_2)$.
Assume now that 
\begin{itemize}
\item $\theta_2 \mapsto f(\theta_1,\theta_2)$ is $\mu$-strongly convex for all $\theta_1$ in~$\Real^{p_1}$;
\item $\theta_1 \mapsto f(\theta_1,\theta_2)$ is differentiable for all~$\theta_2$; 
\item $(\theta_1,\theta_2) \mapsto \nabla_{1} f(\theta_1,\theta_2)$ is $L'$-Lipschitz with respect to~$\theta_1$ and $L$-Lipschitz with respect to~$\theta_2$.\footnote{The notation $\nabla_1$ denotes the gradient with respect to $\theta_1$.}
\end{itemize}
Let us fix $\kappa_1$ in $\Theta_1$. Then, the following function is a majorizing surrogate in $\S_{L''}(\tilde{f},\kappa)$:
\begin{displaymath}
g: \theta_1 \mapsto f(\theta_1,\kappa_2^\star) ~\text{with}~~ \kappa_2^\star \defin \argmin_{\theta_2 \in \Theta_2} f(\kappa_1,\theta_2),
\end{displaymath}
with $L'' = 2L'+ L^{2}/\mu$.  We can indeed apply Lemma~\ref{lemma:danskin},
which ensures that $\tilde{f}$ is differentiable with $\nabla
\tilde{f}(\theta_1) = \nabla_1 f(\theta_1,\theta_2^\star)$ and $\theta_2^\star
\defin \argmin f(\theta_1,\theta_2)$ for all $\theta_1$.  Moreover, $g$ is
$L'$-smooth and $\tilde{f}$ is $L'+L^2/\mu$-smooth according to
Lemma~\ref{lemma:danskin}, and thus $h \defin g-\tilde{f}$ is $L''$-smooth.
Note that a better constant $L''=L'$ can be obtained when~$f$ is convex, as
noted in the appendix of \cite{mairal17}.

The surrogate~$g$ leads to an alternate minimization algorithm; it is then interesting to note that Proposition~\ref{prop:conv2}
provides similar convergence rates as another recent
analysis~\cite{beck2013convergence}, which makes slightly different assumptions
on the function~$f$.
Variational surrogates might also be useful for problems of a single
variable~$\theta_1$. For instance, consider a regression problem with a Huber loss
function~$H$ defined for all $u$ in~$\Real$ as
\begin{displaymath}
   H(u) \defin \left\{ 
      \begin{array}{lr} \frac{u^2}{2\delta} +\frac{\delta}{2} & \text{if}~~  |u| \leq \delta,  \\ 
                      |u|  & \text{otherwise}, \end{array} \right.  \label{eq:huber}
\end{displaymath}
where~$\delta$ is a positive constant.\footnote{To simplify the
notation, we present a shifted version of the traditional Huber loss, which usually 
satisfies $H(0)=0$.} The Huber loss can be seen as a smoothed version of
the~$\ell_1$-norm when~$\delta$ is small, or simply a robust variant of the squared loss
$u \mapsto \frac{1}{2}u^2$ that asymptotically grows linearly.  Then, it is
easy to show that
\begin{displaymath}
   H(u) = \frac{1}{2} \min_{w \geq \delta} \left[\frac{u^2}{w} + w\right]. \label{eq:huber_var}
\end{displaymath} 
Consider now a regression problem with $m$ training data points represented by vectors~$\x_i$ in~$\Real^p$,
associated to real numbers~$y_i$, for $i=1,\ldots,m$. The robust regression problem with the Huber 
loss can be formulated as the minimization over~$\Real^p$ of
\begin{displaymath}
   \tilde{f}: \theta_1 \mapsto  \sum_{i=1}^m H(y_i- \x_i^\top \theta_1) = \min_{\theta_2 \in \Real^m : \theta_2 \geq \delta} \left[f(\theta_1,\theta_2) \defin \frac{1}{2} \sum_{i=1}^m\frac{(y_i- \x_i^\top \theta_1)^2}{\theta_2[i]} +  \theta_2[i]\right],
\end{displaymath}
where $\theta_1$ is the parameter vector of a linear model. The conditions described at the beginning of this section can
be shown to be
satisfied with a Lipschitz constant proportional to~$(1/\delta)$; the resulting
algorithm is the iterative reweighted least-square method, which appears
both in the literature about robust statistics~\cite{lange2}, and about sparse
estimation where the Huber loss is used to approximate the $\ell_1$-norm~\citep{bach8}.

\subsubsection{Jensen surrogates}
Jensen's inequality also provides a natural mechanism to obtain surrogates for convex
functions. Following the presentation of Lange, Hunger and Yang~\citet{lange2}, we consider a convex function
$f: \Real \mapsto \Real$, a vector $\x$ in~$\Real^p$, and define $\tilde{f}:
\Real^p \to \Real$ as $\tilde{f}(\theta) \defin f(\x^\top \theta)$ for
all $\theta$. Let $\w$ be a weight vector in $\Real_+^p$ such that $\|\w\|_1=1$ and $\w[i] \neq 0$ whenever $\x[i] \!\neq\! 0$. Then, we define for any $\kappa$ in~$\Real^p$:
\begin{displaymath}
   g: \theta \mapsto \sum_{i=1}^p \w[i] f \left(\frac{\x[i]}{\w[i]}( \theta[i]-\kappa[i]) + \x^\top \kappa\right),
\end{displaymath}
When $f$ is $L$-smooth, and when $\w[i] \defin |\x[i]|^\nu / \|\x\|_\nu^\nu$, $g$ is in $\S_{L'}(\tilde{f},\kappa)$ with
\begin{itemize}
   \item $L' = L\|\x\|_\infty^2\|\x\|_0$ for $\nu=0$;
   \item $L' = L\|\x\|_\infty\|\x\|_1$ for $\nu=1$;
   \item $L' = L\|\x\|_2^2$ for $\nu=2$.
\end{itemize}
To the best of our knowledge, non-asymptotic convergence rates have not been
studied before for such surrogates, and thus we believe that our analysis may
provide new results in the present case. Jensen surrogates are indeed quite uncommon; they appear
nevertheless in a few occasions. In addition to the few examples given in~\cite{lange2}, they are used for instance in machine
learning by Della Pietra~\citet{pietra} for interpreting boosting procedures
through the concept of \emph{auxiliary functions}.

Jensen's inequality is also used in a different fashion in EM
algorithms~\cite{dempster,neal}. Consider $T$ non-negative functions~$f^t:
\Real^p \mapsto \Real_+$, and, for some~$\kappa$ in~$\Real^p$, define some weights
$\w[t]=  f^t(\kappa) / \sum_{t'=1}^T f^{t'}(\kappa)$. By exploiting the concavity of the logarithm, and assuming hat $\w[t] > 0$ for all~$t$ to simplify, Jensen's
inequality yields
\begin{equation}
   - \log \left( \sum_{t=1}^T f^t(\theta) \right) \leq  - \sum_{t=1}^T \w[t] \log \left(\frac{f^t(\theta)}{\w[t]}\right), \label{eq:em}
\end{equation}
The relation~(\ref{eq:em}) is key to EM algorithms minimizing a negative
log-likelihood. The right side of this equation can be interpreted as a
majorizing surrogate of the left side since it is easy to show that both terms
are equal for~$\theta=\kappa$. Unfortunately the resulting approximation error
functions are not $L$-smooth in general and these surrogates do not follow 
the assumptions of Definition~\ref{def:surrogate_batch}. As a consequence,
our analysis may apply to some EM algorithms, but not to all of them.

\subsubsection{Quadratic surrogates}
When $f$ is twice differentiable and admits a matrix~$\HH$ such that $\HH- \nabla^2 f$ is always positive definite, the following 
function is a first-order majorizing surrogate:
\begin{displaymath}
    g: \theta \mapsto f(\kappa) + \nabla f(\kappa)^\top (\theta-\kappa) + \frac{1}{2}(\theta-\kappa)^\top\HH(\theta-\kappa).
\end{displaymath}
The Lipschitz constant of $\nabla(g-f)$ is the largest eigenvalue of
$\HH- \nabla^2 f(\theta)$ over~$\Theta$.
Such surrogates appear frequently in the statistics and machine learning
literature~\citep{bohning,jebara,khan}. The goal is to to model the global curvature
of the objective function during each iteration, without resorting to the
Newton method. Even though quadratic surrogates do not necessarily lead to
better theoretical convergence rates than simpler Lipschitz gradient
surrogates, they can be quite effective in practice~\cite{jebara}.

\section{An incremental majorization-minimization algorithm: MISO}\label{sec:incremental}
In this section, we introduce an incremental scheme that exploits the
structure~(\ref{eq:prob}) of~$f$ as a large sum of $T$ components.  The most
popular method for dealing
with such a problem when $f$ is smooth and $\Theta=\Real^p$ is probably the
\emph{stochastic gradient descent} algorithm (SGD) and its variants~(see
\cite{nemirovski}). It consists of drawing at iteration~$n$ an index~${\hat
t}_n$ and updating the solution as $\theta_n \leftarrow \theta_{n-1} - \eta_n
\nabla  f^{{\hat t}_n}(\theta_{n-1})$, where the scalar~$\eta_n$ is a
step size.  Another popular algorithm is the \emph{stochastic mirror descent}
algorithm~(see \cite{juditsky}) for general non-smooth convex problems, a
setting we do not consider in this paper since non-smooth functions do not
always admit practical first-order surrogates.

Recently, linear convergence rates for strongly convex functions~$f^t$ have
been obtained in \citet{schmidt2} and \citet{shalev2} by using randomized
incremental algorithms whose cost per iteration is independent of~$T$. The
method SAG \citet{schmidt2} for smooth unconstrained convex optimization is a
randomized variant of the incremental gradient descent algorithm of~Blatt, Hero
and Gauchman~\citep{blatt}, where an estimate of the gradient~$\nabla f$ is incrementally
updated at each iteration. The method SDCA~\citet{shalev2} for strongly
convex composite optimization is a dual coordinate ascent algorithm that
performs incremental updates in the primal~(\ref{eq:prob}).  Unlike SGD, both
SAG and SDCA require storing information about past iterates, which is a key
for obtaining fast convergence rates.

In a different context, incremental EM algorithms have been proposed by Neal
and Hinton~\citet{neal}, where upper bounds of a non-convex negative
log-likelihood function are incrementally updated. By using similar ideas, we introduce
the scheme MISO in Algorithm~\ref{alg:generic_incremental}. At every iteration,
a single function is observed, and an approximate surrogate of~$f$ is
updated. Note that in the same line of work, Ahn et al.~\cite{ahn} have
proposed a block-coordinate descent majorization-minimization algorithm, which
corresponds to MISO when the variational surrogates of
Section~\ref{subsubsec:variational} are used.

\begin{algorithm}[hbtp]
    \caption{Incremental scheme MISO.}\label{alg:generic_incremental}
    \begin{algorithmic}[1]
       \INPUT $\theta_0 \in \Theta$ (initial estimate); $N$ (number of iterations).
    \STATE Initialization: choose some surrogates $g_0^t$ of $f^t$ near $\theta_0$ for all $t$;
    \FOR{ $n=1,\ldots,N$}
    \STATE Randomly pick up one index $\hat{t}_n$ and choose a surrogate $g_n^{\hat{t}_n}$ of $f^{\hat{t}_n}$ near $\theta_{n-1}$; set $g^t_n \defin g^t_{n-1}$ for all $t \neq \hat{t}_n$.
    \STATE Update the solution:
$       \theta_n \in {\displaystyle \argmin_{\theta \in \Theta}} \frac{1}{T} \sum_{t=1}^T g_n^t(\theta)$. 
    \ENDFOR
    \OUTPUT $\theta_{N}$ (final estimate);
    \end{algorithmic}
\end{algorithm}

In the next two sections, we study the convergence properties of the scheme
MISO. We proceed as in Section~\ref{sec:generic}. Specifically, we start with
the non-convex case, focusing on stationary point conditions, and we show that
similar guarantees as for the batch majorization-minimization algorithm hold.
Then, for convex problems, we present convergence rates that essentially apply
to the proximal gradient surrogates. We obtain sublinear rates~$O(T/n)$ for the
general convex case, and linear ones for strongly convex objective functions.
Even though these rates do not show any theoretical advantage over the batch
algorithm, we also present a more surprising result in
Section~\ref{subsec:strong}; in a large sample regime~$T \geq 2L/\mu$, for
$\mu$-strongly convex functions~$f^t$, minorizing surrogates may be used and
faster rates can be achieved.

\subsection{Convergence analysis} \label{subsec:convmiso}
We start our analysis with the non-convex case, and make the following assumption:
\def\assumpBasicB{{\rm (B)}}
\begin{itemize}
   \item[\assumpBasicB] $f$ is bounded below and for all $\theta,\theta'$ in $\Theta$ and all $t$, the
      directional derivative $\nabla f^t(\theta,\theta'-\theta)$ of $f^t$ at~$\theta$ in the direction
      $\theta'-\theta$ exists.
\end{itemize}
Then, we obtain a first convergence result.

\begin{proposition}[\rm Non-convex analysis]\label{prop:conv13}
Assume that~\assumpBasicB~holds and that the surrogates $g_n^{{\hat t}_n}$ from
Algorithm~\ref{alg:generic_incremental} are majorizing~$f^{{\hat t}_n}$ and are in
$\S_{L}(f^{{\hat t}_n},\theta_{n-1})$.
Then, the conclusions of Proposition~\ref{prop:conv1} hold with probability one.
\end{proposition}
\proofatend
We proceed in several steps.

\proofstep{Almost sure convergence of $(f(\theta_n))_{n \geq 0}$}
Let us define $\barg_n \defin \frac{1}{T}\sum_{t=1}^T g_n^t$. We have the
following relation for all $n \geq 1$,
\begin{equation}
   \barg_n = \barg_{n-1} + \frac{g_n^{{\hat t}_n} - g_{n-1}^{{\hat t}_n}}{T}, \label{eq:recur}
\end{equation}
where the surrogates and the index ${\hat t}_n$ are chosen in the algorithm.
Then, we obtain the following inequalities, which hold with probability one for all $n \geq 1$,
\begin{displaymath}
\begin{split}
   \barg_{n}(\theta_n) & \leq \barg_n(\theta_{n-1}) =  \barg_{n-1}(\theta_{n-1}) + \frac{g_n^{{\hat t}_n}(\theta_{n-1}) - g_{n-1}^{{\hat t}_n}(\theta_{n-1})}{T}  \\
                       & = \barg_{n-1}(\theta_{n-1}) + \frac{f^{{\hat t}_n}(\theta_{n-1}) - g_{n-1}^{{\hat t}_n}(\theta_{n-1})}{T} \leq \barg_{n-1}(\theta_{n-1}).
 \end{split}
\end{displaymath}
The first inequality is true by definition of $\theta_n$ and the second one
because $\barg_{n-1}^{{\hat t}_n}$ is a majorizing surrogate of~$f^{{\hat
t}_n}$.  The sequence $(\barg_n(\theta_n))_{ n \geq 0}$ is thus monotonically
decreasing, bounded below with probability one, and thus converges almost
surely. 
By taking the expectation of these previous inequalities, we also obtain that the
sequence $(\E[\barg_{n}(\theta_n)])_{n \geq 0}$ monotonically converges. Thus, the non-positive quantity $\E[f^{{\hat
t}_n}(\theta_{n-1}) - g_{n-1}^{{\hat t}_n}(\theta_{n-1})]$ is the summand
of a converging sum and we have
\begin{displaymath}
\begin{split}
\E\left[ \sum_{n=0}^{+\infty} g_{n}^{{\hat t}_{n+1}}(\theta_{n}) - f^{{\hat t}_{n+1}}(\theta_{n})\right] & = \sum_{n=0}^{+\infty} \E[g_{n}^{{\hat t}_{n+1}}(\theta_{n}) - f^{{\hat t}_{n+1}}(\theta_{n})]  \\
   & = \sum_{n=0}^{+\infty} \E[\E[g_{n}^{{\hat t}_{n+1}}(\theta_{n}) - f^{{\hat t}_{n+1}}(\theta_{n})| \FF_{n}]] \\
   & = \sum_{n=0}^{+\infty} \E[\barg_{n}(\theta_{n}) - f(\theta_{n})] \\
   & =\E\left[\sum_{n=0}^{+\infty} \barg_{n}(\theta_{n}) - f(\theta_{n})\right] < +\infty, \\
\end{split}
\end{displaymath}
where we use Beppo-L\'evy theorem to interchange the expectation and the
sum in front of non-negative quantities, and~$\FF_n$ is the filtration
representing all information up to iteration~$n$ (including~$\theta_n$).
As a result, the sequence $(\barg_{n}(\theta_{n}) - f(\theta_{n}))_{n \geq 0}$ converges almost surely to $0$, implying the almost sure convergence of $(f(\theta_n))_{n \geq 0}$.

\proofstep{Asymptotic stationary point conditions}
Let us define ${\bar h}_n \defin \barg_n-f$, which is $L$-smooth.
Then, for all $\theta$ in $\Theta$ and $n \geq 1$,
\begin{displaymath}
\nabla f(\theta_n,\theta-\theta_n) = \nabla \barg_n(\theta_n,\theta-\theta_n) - \nabla {\bar h}_n(\theta_n)^\top(\theta-\theta_n). 
\end{displaymath}
We have $\nabla \barg_n(\theta_n,\theta-\theta_n) \geq 0$ by definition of $\theta_n$, and $\|\nabla {\bar h}_n(\theta_n)\|_2^2 \leq 2L {\bar h}_n(\theta_n)$, following similar steps as in the proof of Proposition~\ref{prop:conv1}.
Since we have previously shown that $({\bar h}_n(\theta_n))_{n \geq 0}$ almost surely converges to zero, we conclude as in the proof of Proposition~\ref{prop:conv1}, replacing $h_n$ by~$\barh_n$ and~$g_n$ by~$\barg_n$.
\endproofatend

We also give the counterpart of Proposition~\ref{prop:conv1_separable} for Algorithm~\ref{alg:generic_incremental}.
\begin{proposition}[\rm Non-convex analysis - composition]\label{prop:conv13b}
Assume that~\assumpBasicB~is satisfied and that the functions $f^t$ are compositions $f^t = f^{\prime t} \circ e^t$,
   where the functions~$e^t$ are $C$-Lipschitz continuous for some $C>0$.
   Assume also that the functions $g_n^{{\hat t}_n}$ in
   Algorithm~\ref{alg:generic_incremental} are also compositions $g_n^{{\hat t}_n} = g_n^{\prime {\hat t}_n} \circ e^{{\hat t}_n}$, 
   where $g_n^{\prime {\hat t}_n}$ is 
majorizing $f^{\prime {\hat t}_n}$ and is in
$\S_{L}(f^{\prime {\hat t}_n},e^{{\hat t}_n}(\theta_{n-1}))$.
Then, the conclusions of Proposition~\ref{prop:conv13} hold.
\end{proposition}

\proofatend
We first remark that the first part of the proof of
Proposition~\ref{prop:conv13} does not exploit the fact that the approximation
errors $g_n^t - f^t$ are $L$-smooth, but only the fact that $g_n^t$ is
majorizing $f^t$ for all $n$ and $t$. Thus, the first part of the proof of
Proposition~\ref{prop:conv13} holds in the present case, such that $(f(\theta_n))_{n \geq 0}$ almost
surely converges, and the sequence $(\barg_n(\theta_n)- f(\theta_n))_{n \geq 0}$ almost
surely converges to zero, where $\barg_n$ is defined in the proof of
Proposition~\ref{prop:conv13}.

It remains to show that the asymptotic stationary point conditions are
satisfied. To that effect, we follow the proof of
Proposition~\ref{prop:conv1_separable}. We first have, for all $n \geq 1$,
\begin{displaymath}
   \nabla f(\theta_n,\theta-\theta_n) = \nabla \barg_n(\theta_n,\theta-\theta_n) - \frac{1}{T}\sum_{t=1}^T \nabla {\bar h}^t_n(\theta_n,\theta-\theta_n),
\end{displaymath}
with $\nabla \barg_n(\theta_n,\theta-\theta_n) \geq 0$ and $\barh_n^t \defin \barg_n^t-f^t$. Then, following the proof of Proposition~\ref{prop:conv1_separable},
it is easy to show that
\begin{displaymath}
   |\nabla \barh_n^t(\theta_n, \theta-\theta_n)| \leq C\| \nabla \barh_n^{\prime t}(e^t(\theta_n))\|_2{\|\theta-\theta_n\|_2},
\end{displaymath}
where $\barh_n^{\prime t} = \barg_n^{\prime t}-f^{\prime t}$, and we conclude as in Proposition~\ref{prop:conv1_separable}.
\endproofatend

The next lemma provides convergence rates for the convex case, under the
assumption that the surrogate functions are $\rho$-strongly convex with $\rho
\geq L$. The result notably applies to the proximal gradient surrogates of
Section~\ref{subsubsec:proximal}.

\begin{proposition}[\rm Convex analysis for strongly convex surrogate functions]\label{prop:conv16}
Assume that $f$ is convex and bounded below, let $\theta^\star$ be a minimizer of~$f$ on~$\Theta$, and let us define $f^\star \defin \min_{\theta \in \Theta} f(\theta)$. When the surrogates~$g_n^t$ in 
Algorithm~\ref{alg:generic_incremental} are majorizing~$f^t$ and are in $\S_{L,\rho}(f^t,\theta_{n-1})$ with $\rho \geq L$, we have
for all $n \geq 1$,
\begin{equation}
\E[ f(\bartheta_n)-f^\star] \leq \frac{LT\|\theta^\star-\theta_0\|_2^2}{2n}, \label{eq:incr:rate}
\end{equation}
where $\bartheta_n \defin \frac{1}{n} \sum_{i=1}^n \theta_i$ is the average of the iterates.
Assume now that $f$ is $\mu$-strongly convex. For all $n \geq 1$, 
\begin{equation}
   \E[ f(\theta_n) - f^\star] \leq  \left(1- \frac{2\mu}{T(\rho+\mu)}\right)^{n-1}\frac{L\|\theta^\star-\theta_0\|_2^2}{2}.\label{eq:incr:ratemu}
\end{equation}

\end{proposition}
\proofatend
We proceed in several steps. 

\proofstep{Preliminaries}
For all~$n \geq 1$, we introduce the point $\kappa_{n-1}^t$ in $\Theta$ such that $g_n^t$ is in $\S_{L,\rho}(f^t,\kappa_{n-1}^t)$.
We remark that such points are drawn recursively according to the following conditional probability distribution:
$$\PPP(\kappa_{n-1}^t = \theta_{n-1} | \FF_{n-1}) = \delta ~~\text{and}~~ \PPP(\kappa_{n-1}^t = \kappa_{n-2}^t | \FF_{n-1})=1-\delta,$$ 
where $\delta \defin 1/T$, $\FF_n$ is the filtration representing all information up to iteration~$n$ (including~$\theta_n$), and $\kappa_0^t \defin \theta_0$ for all $t$. Thus we have for all $t$ and all $n \geq 1$,
\begin{equation}
   \E[\|\theta^\star-\kappa_{n-1}^t\|_2^2]=   \E[\E[\|\theta^\star-\kappa_{n-1}^t \|_2^2| \FF_{n-1}]] =  \delta\E[\|\theta^\star-\theta_{n-1}\|_2^2] + (1-\delta) \E[\|\theta^\star-\kappa_{n-2}^t\|_2^2].\label{eq:incr:tmp2}
\end{equation}
We also need the following extension of Lemma~\ref{lemma:basic} to the incremental setting: for all~$\theta$ in~$\Theta$ and $n \geq 1$,
\begin{equation}
f(\theta_n) \leq f(\theta) + \frac{1}{T}\sum_{t=1}^T \left( \frac{L}{2}\|\theta-\kappa_{n-1}^t\|_2^2 -\frac{\rho}{2}\|\theta-\theta_n\|_2^2  \right). \label{eq:incr:tmp4}
\end{equation}
The proof of this relation is similar to that of Lemma~\ref{lemma:basic}, exploiting the $\rho$-strong convexity of ${\bar g}_n \defin (1/T)\sum_{t=1}^T g_n^t$.  We
can now study the first part of the proposition.

\proofstep{Non-strongly convex case ($\rho=L$)}
Let us define the quantities $A_n \defin \E[\frac{1}{2T}\sum_{t=1}^T \|\theta^\star-\kappa_n^t\|_2^2]$ and $\xi_n \defin \frac{1}{2}\E[\|\theta^\star-\theta_n\|_2^2]$. Then, we have from~(\ref{eq:incr:tmp4}) with $\theta=\theta^\star$, and by taking the expectation
\begin{displaymath}
    \E[f(\theta_n)-f^\star] \leq LA_{n-1} - L \xi_n.
\end{displaymath}
It follows from~(\ref{eq:incr:tmp2}) that $A_n = \delta\xi_{n} + (1-\delta) A_{n-1}$ and thus, for all $n \geq 1$,
\begin{displaymath}
    \E[f(\theta_n)-f^\star] \leq \frac{L}{\delta}(A_{n-1} - A_n).
\end{displaymath}
By summing the above inequalities, and using Jensen's inequality, we obtain that
\begin{displaymath}
   \E[f(\bartheta_n)-f^\star] \leq  \frac{1}{n}\sum_{i=1}^n\E[f(\theta_i)-f^\star]  \leq \frac{L A_0}{\delta},
\end{displaymath}
leading to the convergence rate of Eq.~(\ref{eq:incr:rate}), since $A_0 = \frac{1}{2}\|\theta^\star-\theta_0\|_2^2$.

\proofstep{$\mu$-strongly convex case}
Assume now that the functions~$f^t$ are $\mu$-strongly convex.
For all~$n \geq 1$, the strong convexity of~$f$ and~(\ref{eq:incr:tmp4}) give us the following inequalities
\begin{displaymath}
  \mu \xi_n \leq \E[f(\theta_n) - f^\star] \leq L A_{n-1} - \rho\xi_n,
\end{displaymath}
Combining this last inequality with~(\ref{eq:incr:tmp2}), we obtain that for all $n \geq 1$,
\begin{displaymath}
   A_n = \delta \xi_n + (1-\delta) A_{n-1} \leq \left( \frac{\delta L}{\mu+\rho} + (1-\delta) \right) A_{n-1}.
\end{displaymath}
Thus, $A_n \leq \beta^n A_0$ with $\beta\defin \frac{ (1-\delta)(\rho+\mu)+\delta L}{\rho+\mu}$.  Since $A_0=\xi_0$, $\E[f(\theta_n) - f^\star] \leq L A_{n-1}$, and $\beta \leq 1- 2\delta \mu/(\rho+\mu)$, we finally have shown the desired convergence rate~(\ref{eq:incr:ratemu}).
\endproofatend

The convergence rate of the previous proposition in the convex case suggests
that the incremental scheme and the batch one of Section~\ref{sec:generic} have
the same overall complexity, assuming that each iteration of the batch
algorithm is $T$ times the one of MISO. For strongly convex functions~$f^t$, we
obtain linear convergence rates, a property shared by SAG or SDCA; it is thus natural to make a
more precise comparison with these other incremental approaches, which we
present in the next two sections. 

\subsection{MISO for smooth unconstrained optimization} \label{subsec:strong}
In this section, we assume that the optimization domain is unbounded---that is,
$\Theta=\Real^p$, and that the functions~$f^t$ are $L$-smooth. When using the
Lipschitz gradient surrogates of Section~\ref{subsubsec:gradient}, MISO amounts
to iteratively using the following update rule:
\begin{equation}
   \theta_n \leftarrow \frac{1}{T} \sum_{t=1}^T \kappa_{n-1}^t -
   \frac{1}{LT}\sum_{t=1}^T \nabla f^t(\kappa_{n-1}^t), \label{eq:MISO1}
\end{equation}
where the vectors $\kappa_{n-1}$ are recursively defined for $n \geq 2$ as
$\kappa_{n-1}^{\hat{t}_{n}}=\theta_{n-1}$ and
$\kappa_{n-1}^{{t}}=\kappa_{n-2}^{{t}}$ for $t \neq \hat{t}_{n}$, with
$\kappa_{0}^t=\theta_0$ for all $t$. It is then easy to see that the complexity
of updating $\theta_n$ is independent of~$T$, by storing the vectors
$\z_n^t=\kappa_{n-1}^t-(1/L)\nabla f^t(\kappa_{n-1}^t)$ and performing the
update $\theta_n = \theta_{n-1} + (1/T)(\z_n^t-\z_{n-1}^t)$.  In comparison,
the approach SAG yields a different, but related, update rule:
\begin{equation}
   \theta_n \leftarrow \theta_{n-1} -
   \frac{\alpha}{T}\sum_{t=1}^T \nabla f^t(\kappa_{n-1}^t), \label{eq:SAG1}
\end{equation}
where the value~$\alpha=1/(16L)$ is suggested in~\cite{schmidt2}. Even
though the rules~(\ref{eq:MISO1}) and~(\ref{eq:SAG1}) seem to be similar to each other at first sight,
they behave differently in practice and do not have the same theoretical
properties. For non-convex problems, MISO is guaranteed to converge, whereas it
is not known whether it is the case for SAG or not.  For convex
problems, both methods have a convergence rate of the same nature---that is, $O(T/n)$.  For
$\mu$-strongly-convex problems, however, the convergence rate of SAG reported
in~\cite{schmidt2} is substantially better than ours. Whereas the expected
objective of SAG decreases with the rate $O(\rho^n)$ with
$\rho_\text{SAG}=1-\min( \mu/(16L), 1/(8T))$, ours decreases with
$\rho_\text{MISO}=1-2\mu/(T(L+\mu))$, which is larger than
$\rho_\text{SAG}$ unless the problem is very well conditioned.

By maximizing the convex dual of~(\ref{eq:prob}) when the functions $f^t$ are
$\mu$-strongly convex, the approach SDCA yields another update rule that
resembles~(\ref{eq:MISO1}) and~(\ref{eq:SAG1}), and offers similar convergence
rates as~SAG. As part of the procedure, SDCA involves large primal gradient
steps~$\theta_{n-1}-(1/\mu)\nabla f^{\hat{t}_n}(\theta_{n-1})$ for updating the
dual variables. It is thus appealing to study whether such large gradient steps
can be used in~(\ref{eq:MISO1}) in the strongly convex case,
regardless of the majorization-minimization principle. In other words, we
want to study the use of the following surrogates within MISO:
\begin{equation}
   g_n^t : \theta \mapsto f^t(\kappa_{n-1}^t) + \nabla f^t(\kappa_{n-1}^t)^\top (\theta-\kappa_{n-1}^t) + \frac{\mu}{2}\|\theta-\kappa_{n-1}^t\|_2^2, \label{eq:lower_surrogates}
\end{equation} 
which are lower bounds of the functions $f^t$ instead of upper bounds. Then,
minimizing $(1/T)\sum_{t=1}^T g_n^t$ amounts to performing the
update~(\ref{eq:MISO1}) when replacing~$L$ by~$\mu$. The resulting algorithm
is slightly different than SDCA, but resembles it. As shown in the next
proposition, the method achieves a fast convergence rate when $T \geq 2L/\mu$,
but may diverge if~$T$ is small. Note that at the same time as us, a similar
result was independently obtained by Defazio et al.~\cite{defazio}, where a
refined analysis provides a slightly better rate, namely the constant $1/3$ in~(\ref{eq:incr:ratemu2}) may be
replaced by~$1/2$.
\begin{proposition}[\rm MISO for strongly-convex unconstrained smooth problems]\label{prop:conv17}
Assume that the functions $f^t$ are $\mu$-strongly convex, $L$-smooth, and bounded below. Let
$\theta^\star$ be a minimizer of~$f$ on~$\Theta$. Assume that $T \geq 2L/\mu$.
When the functions~$g_n^t$ of Eq.~(\ref{eq:lower_surrogates}) are used in
Algorithm~\ref{alg:generic_incremental}, we have for all $n \geq 1$,
\begin{equation}
   \E[ f(\theta_n) - f^\star] \leq \left(1 - \frac{1}{3T}\right)^{n} \frac{2T}{\mu}\|\nabla f(\theta_0)\|_2^2.\label{eq:incr:ratemu2}
\end{equation}
When the functions $f^t$ are lower-bounded by the function~$\theta \mapsto
(\mu/2)\|\theta\|_2^2$, we can use the initialization $\theta_0=0$ and
$g_0^t: \theta \mapsto (\mu/2)\|\theta\|_2^2$ for all $t$. Then, the quantity $({2T}/{\mu})\|\nabla f(\theta_0)\|_2^2$ in~(\ref{eq:incr:ratemu2})
can be replaced by $T f^\star$.
\end{proposition}
\proofatend
As in the proof of
Proposition~\ref{prop:conv13}, we introduce the function $\barg_n \defin
\frac{1}{T}\sum_{t=1}^T g_n^t$, which is
minimized by~$\theta_n$ for $n \geq 1$. Since~$\barg_n$ is a lower bound on~$f$, we have
the relation~$\barg_n(\theta_n) \leq \barg_n(\theta^\star) \leq f^\star$.
Inspired by the convergence proof of SDCA~\cite{shalev2}, which computes
an convergence rate of an expected duality gap, we proceed by studying the
convergence of the sequence $(f^\star-\E[\barg_n(\theta_n)])_{n \geq 1}$.

On the one hand, we have for all $n \geq 1$,
\begin{equation}
   \begin{split}
      \barg_n(\theta_n) & = \barg_n(\theta_{n-1}) - \frac{\mu}{2}\|\theta_n-\theta_{n-1}\|_2^2 \\
                         & = \barg_{n-1}(\theta_{n-1})  + \delta( g_n^{{\hat t}_n}(\theta_{n-1}) - g_{n-1}^{{\hat t}_n}(\theta_{n-1})) - \frac{\mu}{2}\|\theta_n-\theta_{n-1}\|_2^2, 
   \end{split} \label{eq:tmpconv13aa}
\end{equation}
where~$\delta = 1/T$. The first equality is true because $\barg_n$ is quadratic and
is minimized by~$\theta_n$, and the second one uses the relation~(\ref{eq:recur}).
By definition of~$g_n^{{\hat t}}$, we have that $g_n^{{\hat t}_n}(\theta_{n-1}) = f^{{\hat t}_n}(\theta_{n-1})$,
and by taking the expectation, $\E[g_n^{{\hat t}_n}(\theta_{n-1})]=
\E[f^{{\hat t}_n}(\theta_{n-1})] = \E[\E[f^{{\hat t}_n}(\theta_{n-1})|
\FF_{n-1}]] = \E[f(\theta_{n-1})]$, where $\FF_n$ is the 
the filtration representing all information up to iteration~$n$. 
We also have that $\E[g_{n-1}^{{\hat t}_n}(\theta_{n-1})] = \E[\E[g_{n-1}^{{\hat t}_n}(\theta_{n-1}) | \FF_{n-1}]] = \E[\barg_{n-1}(\theta_{n-1})]$.
Thus, we obtain a first useful relation:
\begin{equation}
   \E[\barg_n(\theta_n)] =  (1-\delta)\E[\barg_{n-1}(\theta_{n-1})] + \delta \E[f(\theta_{n-1})] - \frac{\mu}{2}\E\left[\|\theta_n-\theta_{n-1}\|_2^2\right].\label{eq:tmpconv13a}
\end{equation}
On the other hand, for all $n\geq 2$, 
\begin{equation}
   \begin{split}
      \barg_n(\theta_n) & = \barg_{n-1}(\theta_n) + \delta (g_n^{{\hat t}_n}(\theta_{n}) - g_{n-1}^{{\hat t}_n}(\theta_{n})) \\
                        & = \barg_{n-1}(\theta_{n-1}) \!+\! \frac{\mu-\delta L}{2}\|\theta_{n}\!-\!\theta_{n-1}\|_2^2 \!+\! \delta \left(g_n^{{\hat t}_n}(\theta_{n}) \!+\! \frac{L}{2}\|\theta_n\!-\!\theta_{n-1}\|_2^2 - g_{n-1}^{{\hat t}_n}(\theta_{n})\right) \\
                       & \geq \barg_{n-1}(\theta_{n-1}) + \frac{\mu-\delta L}{2}\|\theta_{n} - \theta_{n-1}\|_2^2.
   \end{split}  \label{eq:tmpconv13c}
\end{equation}
We have used the fact that $\theta \mapsto g_n^{{\hat t}_n}(\theta) +
(L/2)\|\theta-\theta_{n-1}\|_2^2$ is a majorizing surrogate of $f^{{\hat
t}_n}$, whereas $g_{n-1}^{{\hat t}_n}$ is minorizing~$f^{{\hat t}_n}$.
By adding twice~(\ref{eq:tmpconv13c}) after taking the expectation
and once~(\ref{eq:tmpconv13a}), 
we obtain that for all $n \geq 2$,
\begin{equation}
   \begin{split}
      3\E[\barg_n(\theta_n)] & \geq (3 - \delta)\E[\barg_{n-1}(\theta_{n-1})] + \delta\E[ f(\theta_{n-1})] + \left(\frac{\mu}{2}-\delta L\right)\E[\|\theta_n-\theta_{n-1}\|_2^2]  \\
         & \geq  (3 - \delta)\E[\barg_{n-1}(\theta_{n-1})] + \delta\E[ f(\theta_{n-1})],
   \end{split}
      \label{eq:tmpconv13b}
\end{equation}
where the second inequality comes from the large sample size condition $\delta L \leq \mu/2$.
Since $\E[f(\theta_{n-1})] \geq f^\star$, this immediately gives for $n\geq 2$,
\begin{displaymath}
f^\star - \E\left[\barg_n(\theta_n)\right] \leq \left(  1- \frac{1}{3T} \right) \left(f^\star - \E\left[\barg_{n-1}(\theta_{n-1})\right]\right).
\end{displaymath}
To obtain a convergence rate for $\E[f(\theta_n)]-f^\star$, we use again Eq.~(\ref{eq:tmpconv13b}). For $n \geq 2$,
\begin{equation}
   \begin{split}
      \delta(\E[ f(\theta_{n-1})] - f^\star) & \leq \delta(\E[ f(\theta_{n-1})] - \E[\barg_{n-1}(\theta_{n-1})]) \\ 
        & \leq 3 (\E[\barg_n(\theta_n)]- \E[\barg_{n-1}(\theta_{n-1})]) \\
        & \leq 3 (f^\star- \E[\barg_{n-1}(\theta_{n-1})]) \\
        & \leq 3 \left(  1- \frac{1}{3T} \right)^{n-2}\left(f^\star - \barg_{1}(\theta_{1})\right),
   \end{split} \label{eq:tmpincr1}
\end{equation}
and we obtain the convergence rate (\ref{eq:incr:ratemu2}) by first noticing that
\begin{displaymath}
   \begin{split}
      f^\star - \barg_1(\theta_1) & = f^\star - \barg_1(\theta_0) + \frac{\mu}{2}\|\theta_0-\theta_1\|_2^2 \\
                                  & = f^\star - f(\theta_0) + \frac{\mu}{2}\left\|\frac{1}{\mu} \nabla f(\theta_0)\right\|_2^2 \\
      & \leq \frac{1}{2\mu}\|\nabla f(\theta_0)\|_2^2,
   \end{split}
\end{displaymath}
where we use the relation $\barg_1 = \barg_0$
and~$\barg_0(\theta_0)=f(\theta_0)$. Then, we use the fact that $(1-1/3T) \geq
5/6$ since $T \geq 2L/\mu \geq 2$, such that $3(1-1/3T)^{-1}/(2\mu)  \leq
9/(5\mu) \leq 2/\mu$.

To prove the last part of the proposition, we remark that all inequalities we
have proved so far for $n\geq 2$, become true for $n=1$. Thus, the
last inequality in~(\ref{eq:tmpincr1}) is also true when replacing $n-2$ by
$n-1$ and $\barg_1(\theta_1)$ by~$\barg_0(\theta_0)=0$. 
\endproofatend

The proof technique is inspired in part by the one of SDCA~\cite{shalev2}; the
quantity $\sum_{t=1}^T g_n^t(\theta_n)$ is indeed a lower bound of $f^\star$,
and plays a similar role as the dual value in SDCA. We remark that the
convergence rate~(\ref{eq:incr:ratemu2}) improves significantly upon the 
original one~(\ref{eq:incr:ratemu}), and is similar to the one of SAG when $T$
is larger than $2L/\mu$.\footnote{Note that a similar assumption appears in the
first analysis of SAG published in~\cite{leroux} before its refinement
in~\cite{schmidt2}.} However, Proposition~\ref{prop:conv17} only applies to
strongly convex problems. In other cases, the more conservative
rule~(\ref{eq:MISO1}) should be preferred in theory, even though we present
heuristics in Section~\ref{subsec:heuristics} that suggest using larger
step sizes than $1/L$ in practice.

\subsection{MISO for composite optimization}
When $f$ can be written as $f=(1/T)\sum_{t=1}^T f_1^t +  f_2$, where the
functions $f_1^t$ are $L$-smooth, we can use the proximal gradient surrogate
presented in Section~\ref{subsubsec:proximal}; it yields the following rule:
\begin{equation}
   \theta_n \in \argmin_{\theta \in \Theta}  \frac{1}{2}\left\| \theta -
   \left(\frac{1}{T}\sum_{t=1}^T\kappa_{n-1}^t - \frac{1}{LT} \sum_{t=1}^T
   \nabla f_1^t ( \kappa_{n-1}^t)\right)  \right\|_2^2 +
   \frac{\lambda}{L}f_2(\theta), \label{eq:MISO_composite}
\end{equation}
where the vectors~$\kappa_{n-1}^t$ are defined as in
Section~\ref{subsec:strong}. This update is related to SDCA, as well as to
stochastic methods for composite convex optimization such as the
regularized dual averaging algorithm of Xiao~\cite{xiao}. As in the previous
section, we obtain guarantees for non-convex optimization, but our linear
convergence rate for strongly convex problems is not as fast as the one of
SDCA. Even though we do not have a similar result as
Proposition~\ref{prop:conv17} for the composite setting, we have observed that
using a smaller value for~$L$ than the theoretical one could work well in
practice. We detail such an empirical strategy in the next section.

\subsection{Practical implementation and heuristics}\label{subsec:heuristics}
We have found the following strategies to improve the practical performance
of MISO.

\paragraph{Initialization} A first question is how to
initialize the surrogates~$g_0^t$ in practice. Even though
we have suggested the functions $g_0^t$ to be in~$\S_{L}(f^t,\theta_0)$ in
Algorithm~\ref{alg:generic_incremental}, our analysis weakly
relies on this assumption. In fact, most of our results hold when choosing 
surrogates computed at points~$\kappa_0^t$ that are not necessarily equal to~$\theta_0$; at most
only constants from the convergence rates would be affected by such a change.  An 
effective empirical strategy is inspired by the second part of
Proposition~\ref{prop:conv17}: we first define functions~$g_{0}^t: \theta
\mapsto (L/2)\|\theta-\theta_0\|_2^2$, and perform $T$ iterations of MISO without randomization,
selecting the function~$f^t$ at iteration~$t$, such that
each surrogate is updated exactly once. Then, we use these updated surrogates
for initializing the regular randomized scheme.

\paragraph{Warm restart and continuation}
When available, warm restart can be used for initializing the surrogates.
Assume that we are interested in minimizing a composite function
$(1/T)\sum_{t=1}^T f_1^t(\theta) + \lambda f_2(\theta)$, which is parameterized
by a scalar~$\lambda$, and that we want to obtain a minimizer for several
parameter values $\lambda_1 < \lambda_2 < ... < \lambda_M$. We first solve the
problem for $\lambda=\lambda_M$, and then use the surrogates obtained at the
end of the optimization for initializing the algorithm when addressing the
problem with $\lambda=\lambda_{M-1}$. We proceed similarly going from larger to
smaller values of $\lambda$. We have empirically observed that the warm restart
strategy could be extremely efficient in practice, and would deserve further
study in a future work.

\paragraph{Heuristics for selecting step sizes}
Choosing proximal gradient surrogates~$g^t$ requires choosing some Lipschitz
constant~$L$ (or a strong convexity parameter $\mu$ for
Proposition~\ref{prop:conv17}), which leads to a specific step size
in~(\ref{eq:MISO_composite}). However, finding an appropriate
step size can be difficult in practice for several reasons:
(i) in some cases, these parameters are unknown; (ii) even though a global Lipschitz
constant might be available, a local Lipschitz constant could be more effective; 
(iii) the convergence rates of Proposition~\ref{prop:conv16} can be obtained by
choosing a smaller value for~$L$ than the ``true'' Lipschitz constant, as long
as the inequality $\E[f(\theta_n)] \leq \E[\barg_n(\theta_n)]$ is always
satisfied, where $\barg_n \defin (1/T)\sum_{t=1}^T g_n^t$.
This motivates the following heuristics:
\begin{itemize}
   \item[MISO1] first perform one pass over $\eta\!=\!5\%$ of the data to
      select a constant~$L_1 = 2^{-k}L_0$ with~$k$ chosen among
      positive integers, yielding the smallest objective on the data subset, where~$L_0$ 
      is an upper bound of the true Lipschitz constant.
   \item[MISO2] proceed as in MISO1, but choose a more aggressive strategy~$L_2=L_1 \eta$;
      during the optimization, 
      compute the quantities $a_n^t$ and $b_n^t$ defined as
      $a_n^t = a_{n-1}^t$, $b_n^t = b_{n-1}^t$ if $t \neq \hat{t}_n$, and
      otherwise $a_n^{{\hat t}_n} = f^{{\hat t}_n}(\theta_{n-1})$, $b_n^{{\hat
      t}_n} = g_{L_2}^{{\hat t}_n}(\theta_{n-1})$, where we have parameterized the surrogates $g^t$ by $L_2$.
      Every $T$ iterations, compare the sums $A_n = \sum_{t=1}^T a_n^t$ and $B_n = \sum_{t=1}^T b_n^t$. 
      If $A_n \leq B_n$, do nothing; otherwise, increase the value of $L_2$ until this inequality
      is satisfied.
\end{itemize}
The heuristic MISO2 is more aggressive than MISO1 since it starts with a
smaller value for $L$. After every iteration, this value is possibly increased
such that on average, the surrogates ``behave'' as majorizing functions.
Even though this heuristic does not come with any theoretical guarantee, it was
found to perform slightly better than MISO1 for strongly-convex problems.

\paragraph{Using a different parameter $L_t$ for every function $f_t$}
Even though our analysis was conducted with a global parameter~$L$ for
simplicity, it is easy to extend the analysis when the parameter $L$ is
adjusted individually for every surrogate. This is useful when the functions
$f_t$ are heterogeneous.

\paragraph{Parallelization with mini-batches}
The complexity of MISO is often dominated by the cost of updating the
surrogates $g_n^{{\hat t}_n}$, which typically requires computing the gradient
of a function.  A simple extension is to update several surrogates at the same
time, when parallel computing facilities are available.

\section{Experimental validation}\label{sec:exp}
In this section, we evaluate MISO on large-scale machine learning problems.
Our implementation is coded in C++ interfaced with Matlab and is freely 
available in the open-source software package
SPAMS~\cite{mairal7}.\footnote{\url{http://spams-devel.gforge.inria.fr/}.} All
experiments were conducted on a single core of a 2GHz Intel CPU with $64$GB of~RAM.

\paragraph{Datasets}
We use six publicly available datasets, which consist of
pairs~$(y_t,\x_t)_{t=1}^T$, where the $y_t$'s are labels in~$\{-1,+1\}$
and the $\x_t$'s are vectors in~$\Real^p$ representing data points.
The datasets are described in Table~\ref{table:datasets}. 
\textsf{alpha}, \textsf{rcv1}, \textsf{ocr}, and \textsf{webspam} are
obtained from the 2008 Pascal large-scale learning
challenge.\footnote{\url{http://largescale.ml.tu-berlin.de}.}
\textsf{covtype} and \textsf{real-sim} are obtained from the LIBSVM
website.\footnote{\url{http://www.csie.ntu.edu.tw/~cjlin/libsvm/}.} 
The datasets are pre-processed as follows: all dense datasets 
are standardized to have zero-mean and unit variance for every
feature. The sparse datasets are normalized such that each $\x_t$
has unit $\ell_2$-norm.

\begin{table}
   \centering
\caption{Description of datasets used in our experiments.}
\label{table:datasets}
\begin{tabular}{|l|c|c|c|c|c|}
\hline
name & $T$ & $p$ & storage & density & size (GB) \\ 
\hline
\textsf{covtype} & $581\,012$ & $54$ & dense & 1 & $0.23$ \\
\hline
\textsf{alpha} & $500\,000$ & $500$ & dense & 1 & $1.86$\\
\hline
\textsf{ocr} & $ 2\,500\,000$ & $1\,155$ & dense & 1 & $21.5$ \\
\hline
\textsf{real-sim} & $ 72\,309$ & $20\,958$ & sparse & 0.0024 & $0.056$ \\
\hline
\textsf{rcv1} & $781\,265$ & $47\,152$ & sparse & 0.0016 & $0.89$ \\
 \hline
 \textsf{webspam} & $250\,000$  & $16\,091\,143$ & sparse & 0.0002  & $13.90$  \\
\hline
\end{tabular}
\end{table}

\subsection{$\ell_2$-logistic regression}\label{subsec:explog}
We consider the $\ell_2$-regularized 
logistic regression problem, which can be formulated as follows:
\begin{equation}
\min_{\theta \in \Real^p} \frac{1}{T} \sum_{t=1}^T \ell(y_t ,\x_t^{\top}\theta) + \frac{\lambda}{2}\|\theta\|_2^2, \label{eq:logistic}
\end{equation}
where~$\ell(u,\hat{u}) = \log(1+e^{-u \hat{u}})$ for all $(u,\hat{u})$.
Following~\cite{schmidt2}, we report some results obtained with different
methods with the parameter $\lambda=1/T$, which is argued to be of the same
order of magnitude as the smallest value that would be used in practice for
machine learning problems. We also performed experiments with the
values~$\lambda=0.1/T$ and~$\lambda=10/T$ to study the impact of the strong
convexity parameter; the output of these two additional experiments is not
reported in the present paper for space limitation reasons, but it will
be discussed and taken into account in our conclusions. 
The algorithms included in the comparison are: 
\begin{itemize}
   \item[SGD-h] the stochastic gradient descent algorithm with a heuristic for
      choosing the step-size similar to MISO1, and inspired by Leon Bottou's sgd toolbox
      for machine learning.\footnote{available here: \url{http://leon.bottou.org/projects/sgd}.}
      A step-size of the form $\rho/\sqrt{n+n_0}$ is automatically adjusted when performing one
      pass on $\eta=5\%$ of the training data. We obtain consistent results with the performance 
      of SGD reported by Schmidt et al.~\cite{schmidt2} when the step-size is chosen from hindsight. 
      Based on their findings, we do not include in our figures other variants
      of SGD, {\it e.g.}, \cite{duchi4,ghadimi2,hazan2,xiao}.
   \item[FISTA] the accelerated gradient method proposed by Beck and Teboulle~\cite{beck} with
      a line-search for automatically adjusting the Lipschitz constant.
   \item[SDCA] the algorithm of Shalev-Schwartz and Zhang~\cite{shalev2},
      efficiently implemented in the language C by Mark Schmidt.\footnote{available here:
      \url{http://www.di.ens.fr/~mschmidt/Software/SAG.html}.}
   \item[SAG] a fast implementation in C also provided by Mark Schmidt~\cite{schmidt2}. We use  
      the step-size $1/L$ since it performed similar to their heuristic line search.
   \item[MISO0] the majorization-minimization algorithm MISO, using the trivial upper bound $L^t=0.25\|\x_t\|_2^2$ on the Lipschitz
      constant for example~$t$.
   \item[MISO1] the majorization-minimization heuristic MISO1 described in Section~\ref{subsec:heuristics}.
   \item[MISO2] the heuristic MISO2, also described in Section~\ref{subsec:heuristics}.
   \item[MISO$\mu$] the update rule corresponding to Proposition~\ref{prop:conv17}.
\end{itemize}
For sparse datasets, MISO0, MISO1, and MISO2 are not practical since
they suffer from a $O(Tp)$ memory cost. Their update rules
can indeed be rewritten 
\begin{equation*}
   \theta_n \leftarrow \theta_{n-1} - \frac{1}{T}\left( \left(\theta_{n-1} - \frac{1}{L}\nabla f^{{\hat t}_n}(\theta_{n-1})\right) - \left( \kappa_{n-1}^{{\hat t}_n}  - \frac{1}{L}\nabla f^{{\hat t}_n}(\kappa_{n-1}^{{\hat t}_n}) \right)\right), 
\end{equation*}
where $f^t: \theta \mapsto  \ell(y_t ,\x_t^{\top}\theta) +
\frac{\lambda}{2}\|\theta\|_2^2$. Thus, for every
example~$t$, the algorithm requires storing  the dense vector $\kappa_{n-1}^{t}  - ({1}/{L})\nabla
f^{t}(\kappa_{n-1}^{{t}})$.  Therefore, we use mini-batches of size
$\lfloor1/d\rfloor$, where $d$ is the density of the dataset; the resulting
algorithms, which we denote by~MISO0-mb, MISO1-mb, and MISO2-mb, have a storage cost equal
to $O(d pT)$, which is the same as the dataset.

On the other hand, the update rule MISO$\mu$ applied to the $\lambda$-strongly convex functions $f^t$ admits a simpler and computationally~cheaper~form.
Since~$\kappa_{n-1}^{t}  - ({1}/{\lambda})\nabla f^{t}(\kappa_{n-1}^{{t}}) = -(1/\lambda)\ell'(y_{{t}} ,\x_{{t}}^{\top}\kappa_{n-1}^{{t}})  \x_{{t}}$, the update becomes
\begin{equation}
   \theta_n \leftarrow \theta_{n-1} - \frac{1}{T\lambda}\left( \ell'(y_{\hat{t}_n} ,\x_{\hat{t}_n}^{\top}\theta_{n-1}) - \ell'(y_{\hat{t}_n} ,\x_{\hat{t}_n}^{\top}\kappa_{n-1}^{{\hat t}_n})   \right) \x_{\hat{t}_n}, \label{eq:MISO4b}
\end{equation}
where $\ell'$ denotes the derivative of~$\ell$ with respect to its second
argument. Assuming that the dataset fits into memory, the only extra quantities
to store are the scalars~$\ell'(y_{\hat{t}_n}
,\x_{\hat{t}_n}^{\top}\kappa_{n-1}^{{\hat t}_n})$, and the resulting memory
cost is simply $O(T)$.

\begin{figure}[hbtp]
   \centering
   \includegraphics[width=0.33\linewidth]{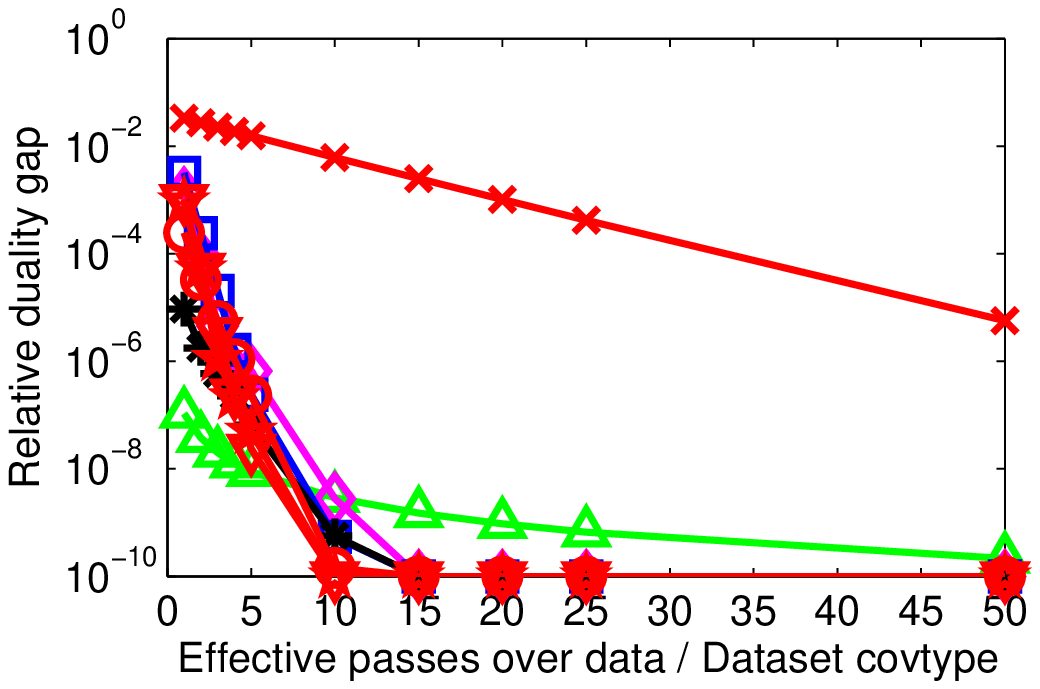}\hfill
   \includegraphics[width=0.33\linewidth]{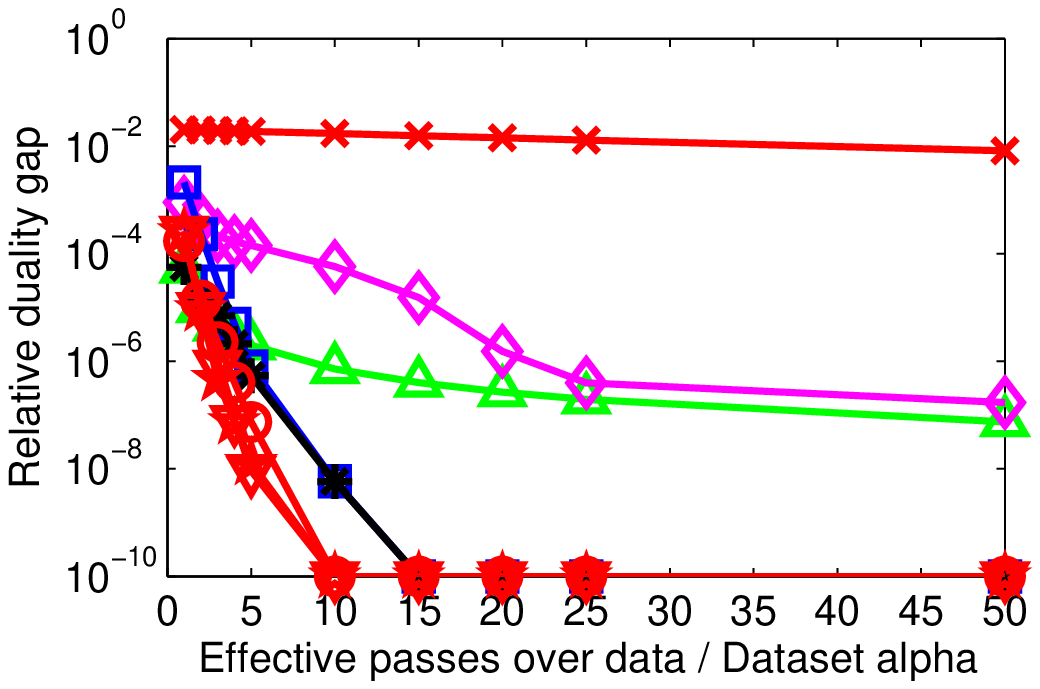}\hfill
   \includegraphics[width=0.33\linewidth]{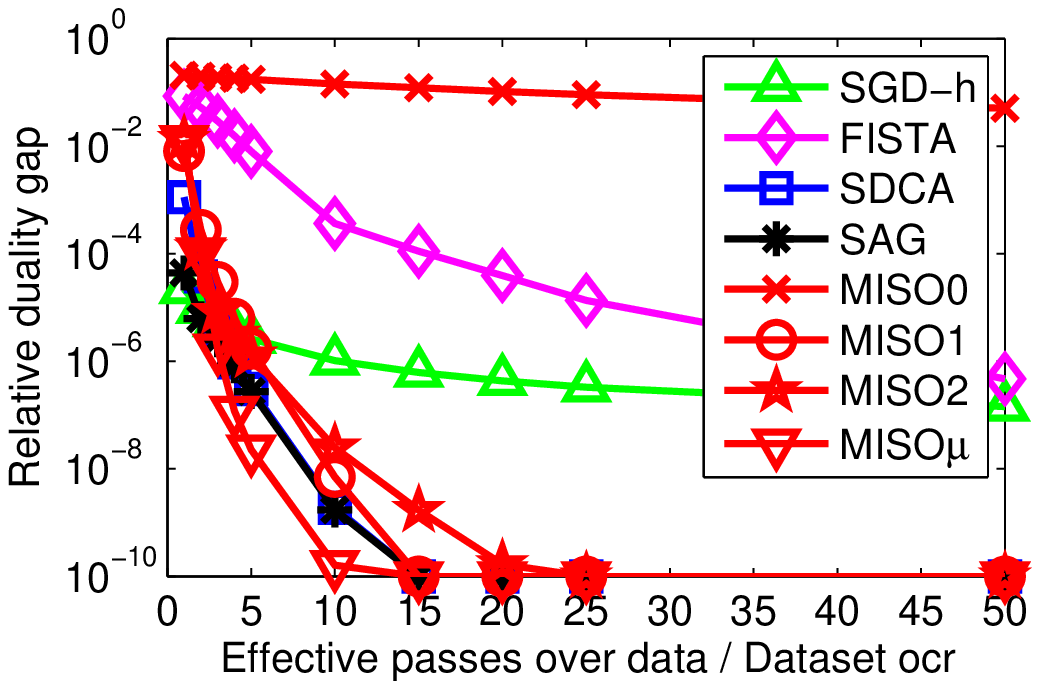}\\
   \includegraphics[width=0.33\linewidth]{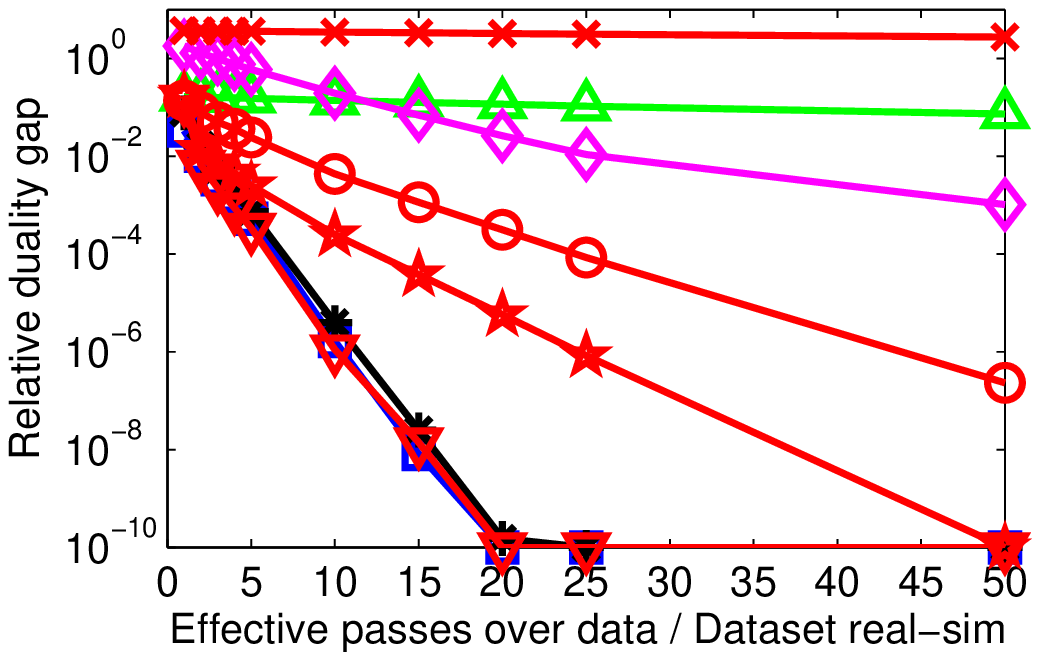}\hfill
   \includegraphics[width=0.33\linewidth]{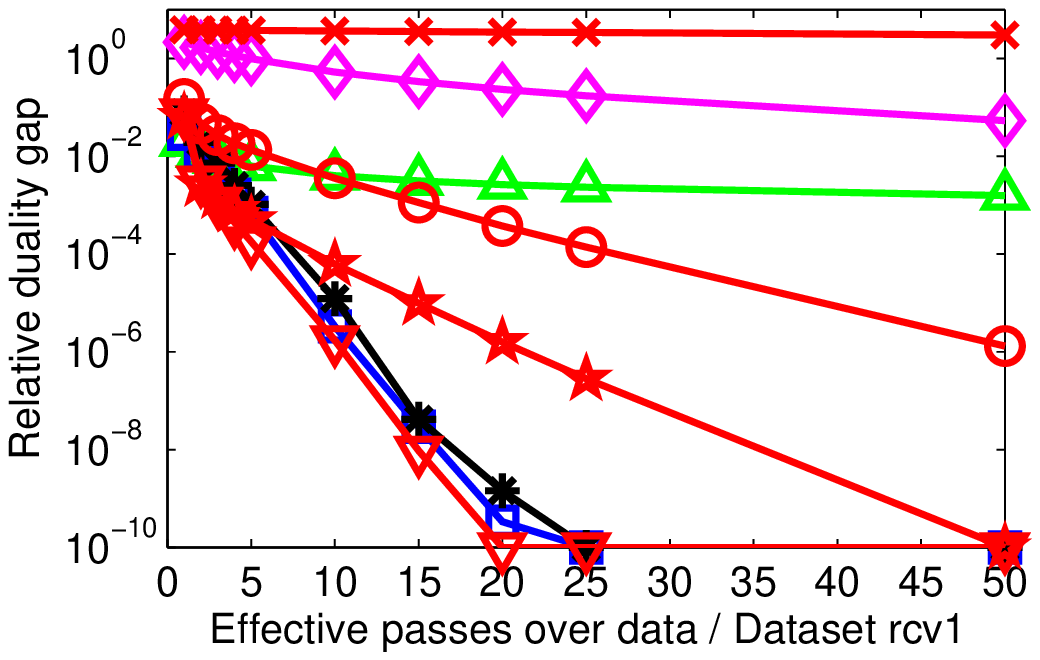}\hfill
   \includegraphics[width=0.33\linewidth]{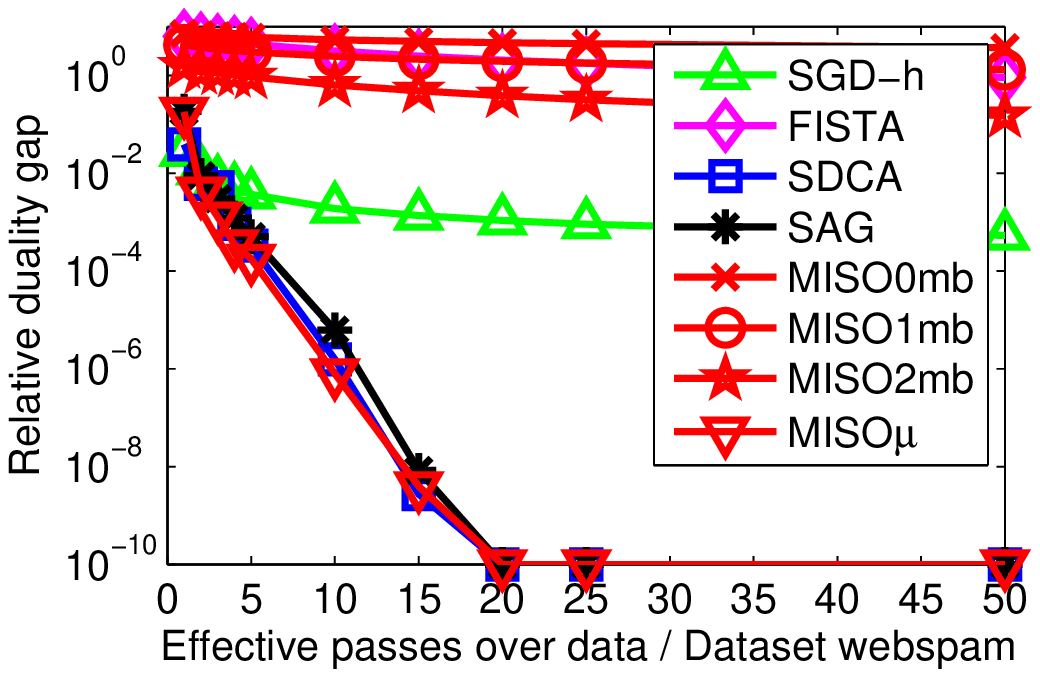}
   \caption{Relative duality gap obtained for logistic regression with respect to the number of passes over the data.}
   \label{fig:l2epochs}
\end{figure}

\begin{figure}[hbtp]
   \centering
   \includegraphics[width=0.33\linewidth]{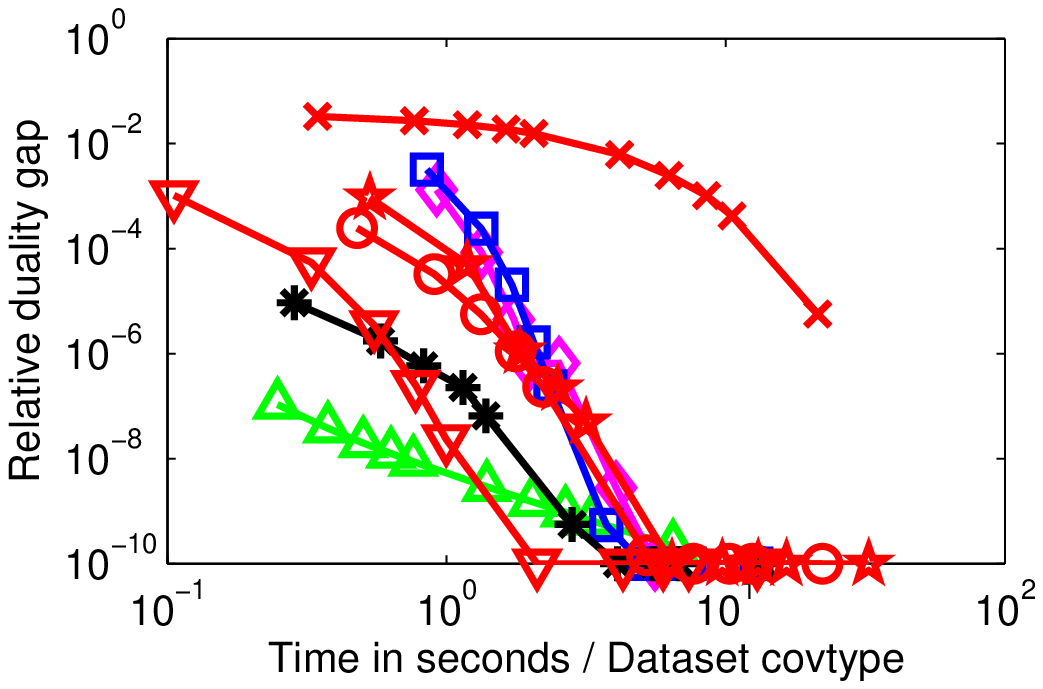}\hfill
   \includegraphics[width=0.33\linewidth]{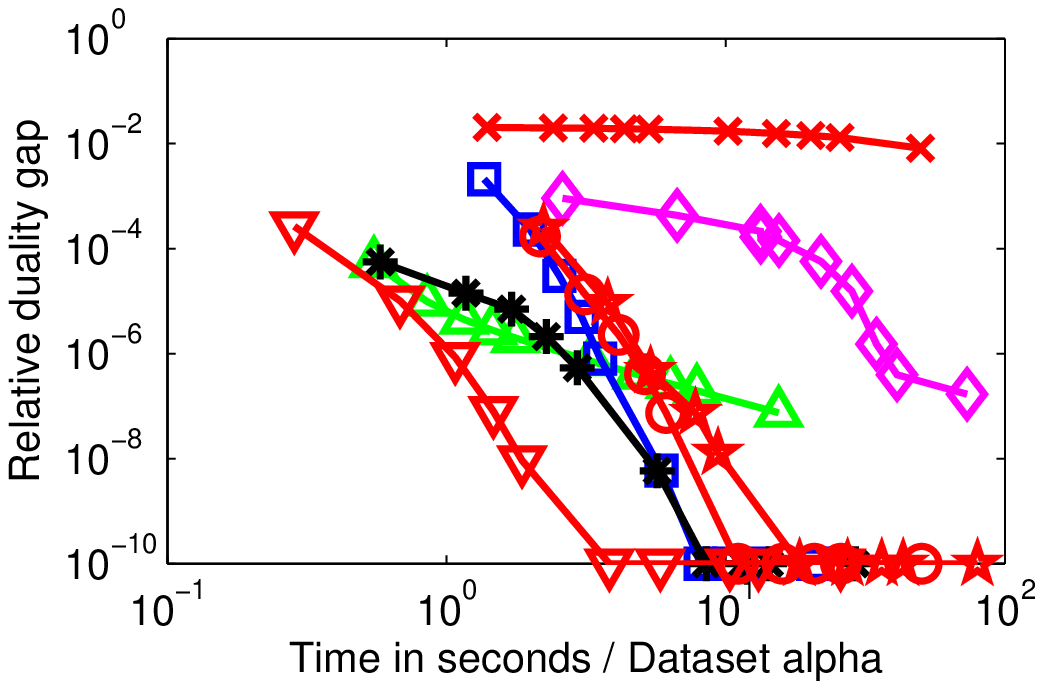}\hfill
   \includegraphics[width=0.33\linewidth]{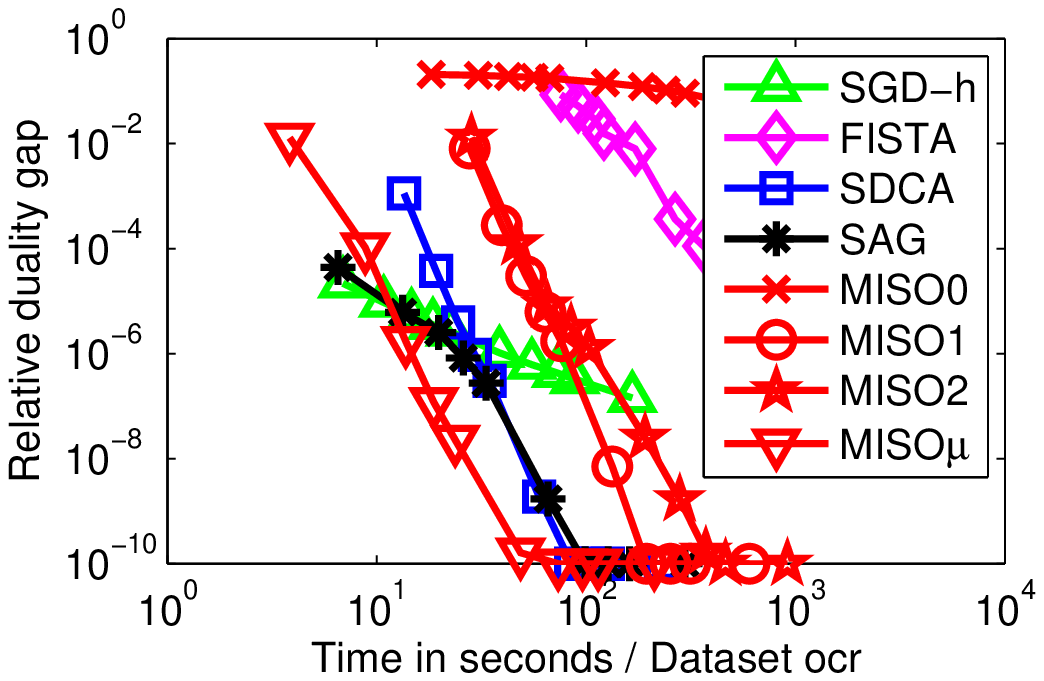}\\
   \includegraphics[width=0.33\linewidth]{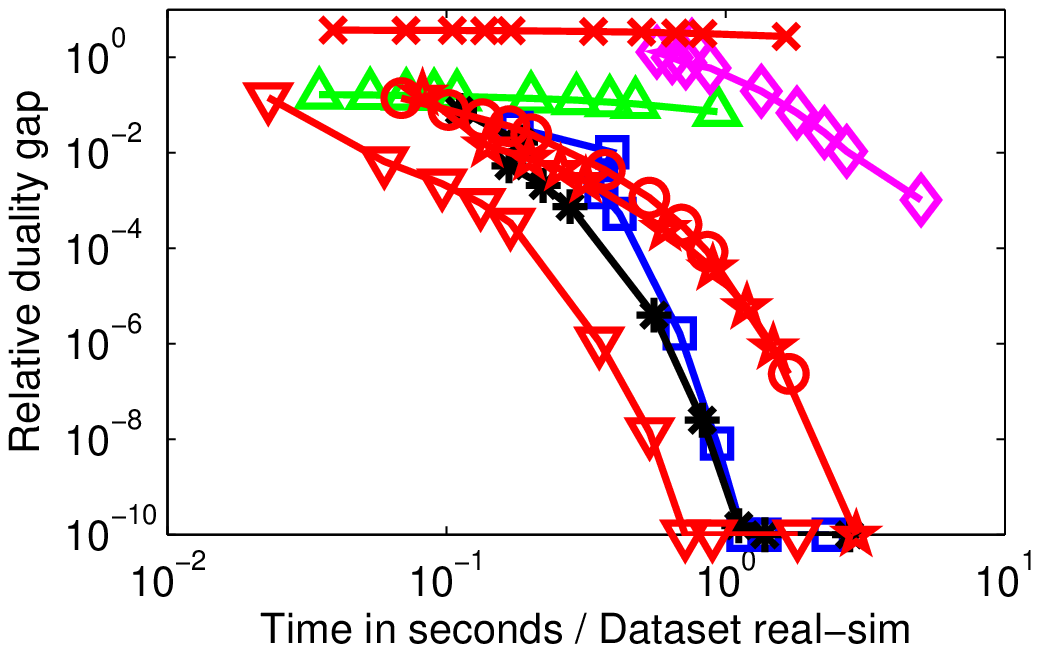}\hfill
   \includegraphics[width=0.33\linewidth]{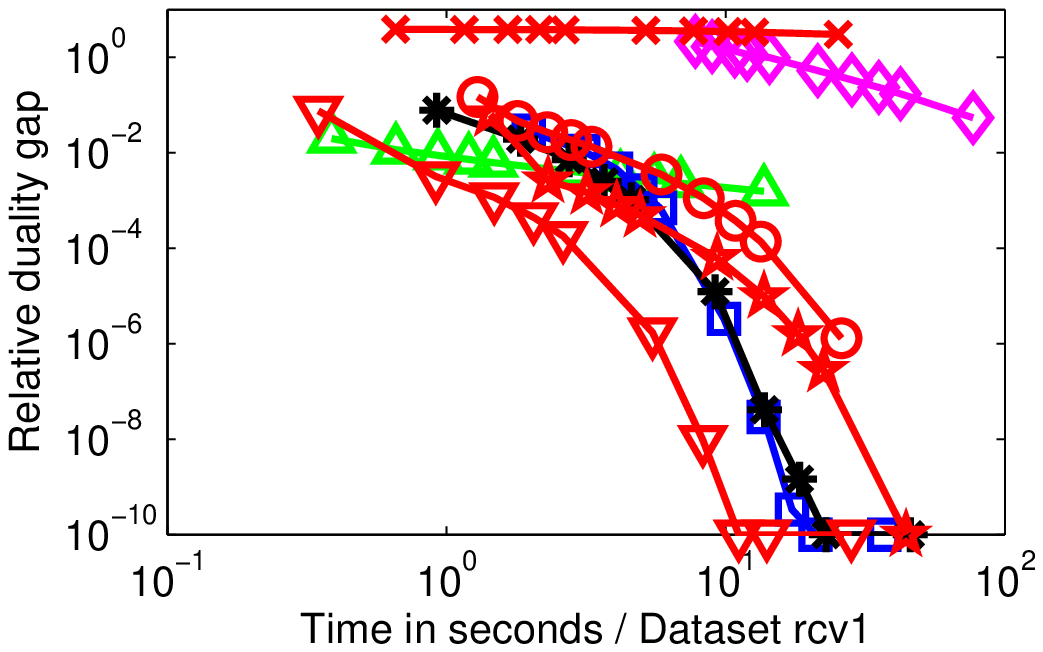}\hfill
   \includegraphics[width=0.33\linewidth]{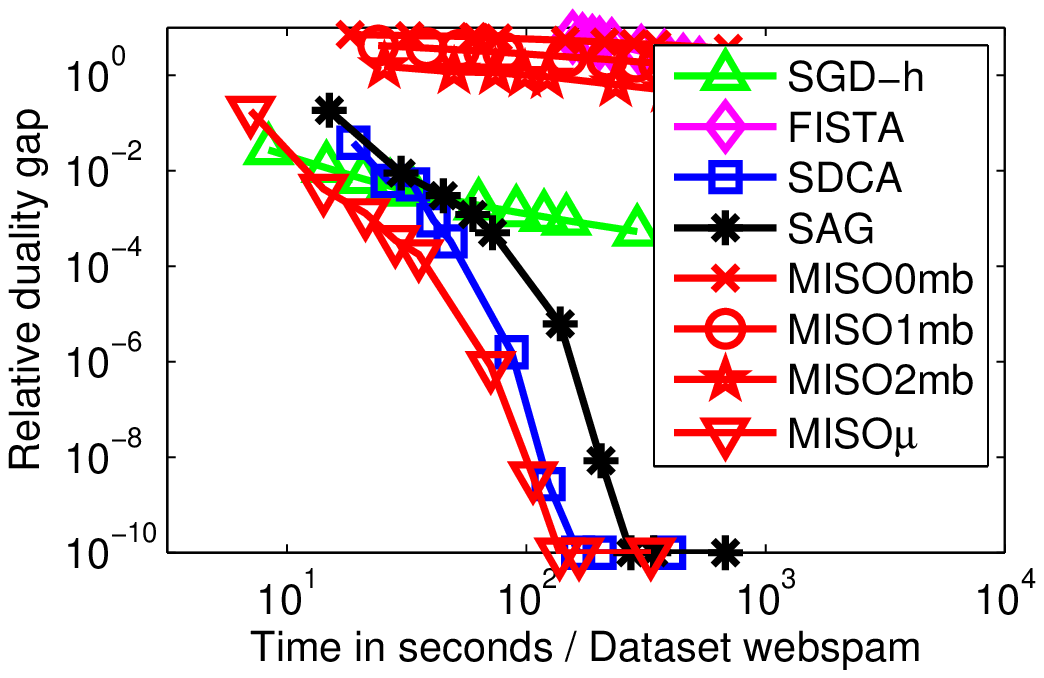}
   \caption{Relative duality gap obtained for logistic regression
   with respect to the CPU time.} \label{fig:l2time}
\end{figure}

We present our comparison of the above methods with~$\lambda=1/T$ on
Figures~\ref{fig:l2epochs} and~\ref{fig:l2time}, where we plot the relative
duality gap defined as $(f(\theta_n)-g^\star)/g^\star$, where $g^\star$ is the
best value of the Fenchel dual that we have obtained during our experiments.
The conclusions of our empirical study are the following:
\begin{itemize}
   \item {\it SAG, SDCA and MISO$\mu$:} these methods perform similarly and were consistently 
      the fastest, except in the regime $T < 2L/\mu$ where MISO$\mu$ can diverge;
   \item {\it the four variants of MISO:} as predicted by its theoretical
      convergence rate, MISO0 does not perform better than ISTA~\cite{beck} without line-search (not reported in the figures). MISO1
      and MISO2 perform significantly better. MISO$\mu$ is always 
      better or as good as MISO1 and MISO2, except for sparse datasets with
      $\lambda=0.1/T$ where the condition~$T \geq
      2L/\mu$ is not satisfied; 
   \item {\it influence of mini-batch:} whereas MISO2 performs equally well as SAG/SDCA for dense datasets, 
      mini-batches for sparse datasets makes it slower;
   \item {\it stochastic gradient descent:} SGD-h performs always well at the beginning
      of the procedure, but is not competitive compared to incremental approaches after a few
      passes over the data.
\end{itemize}
Note that an evaluation of a preliminary version of MISO2 is presented
in~\cite{mairal17} for the $\ell_1$-regularized logistic regression problem,
where the objective function is not strongly convex. Our experimental findings
showed that MISO2 was competitive with state-of-the-art solvers based on
active-set and coordinate descent algorithms~\cite{fan2}.

\subsection{Non-convex sparse estimation}
\begin{figure}[hbtp]
   \centering
   \includegraphics[width=0.33\linewidth]{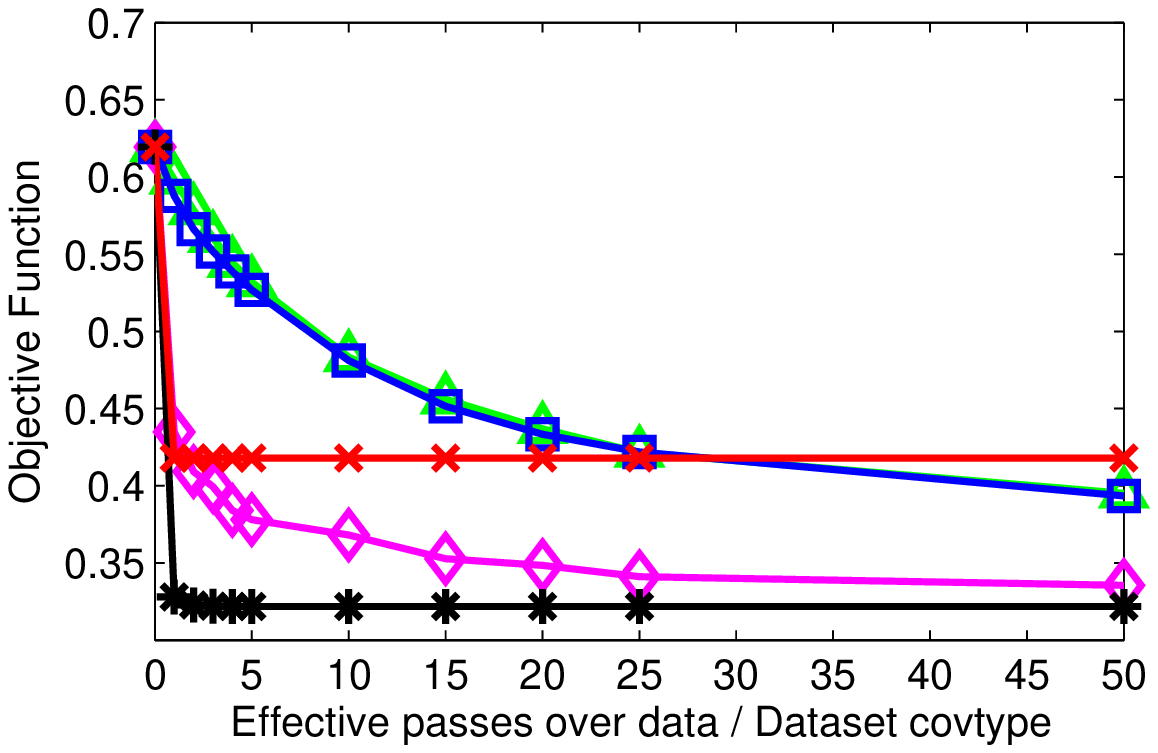}\hfill
   \includegraphics[width=0.33\linewidth]{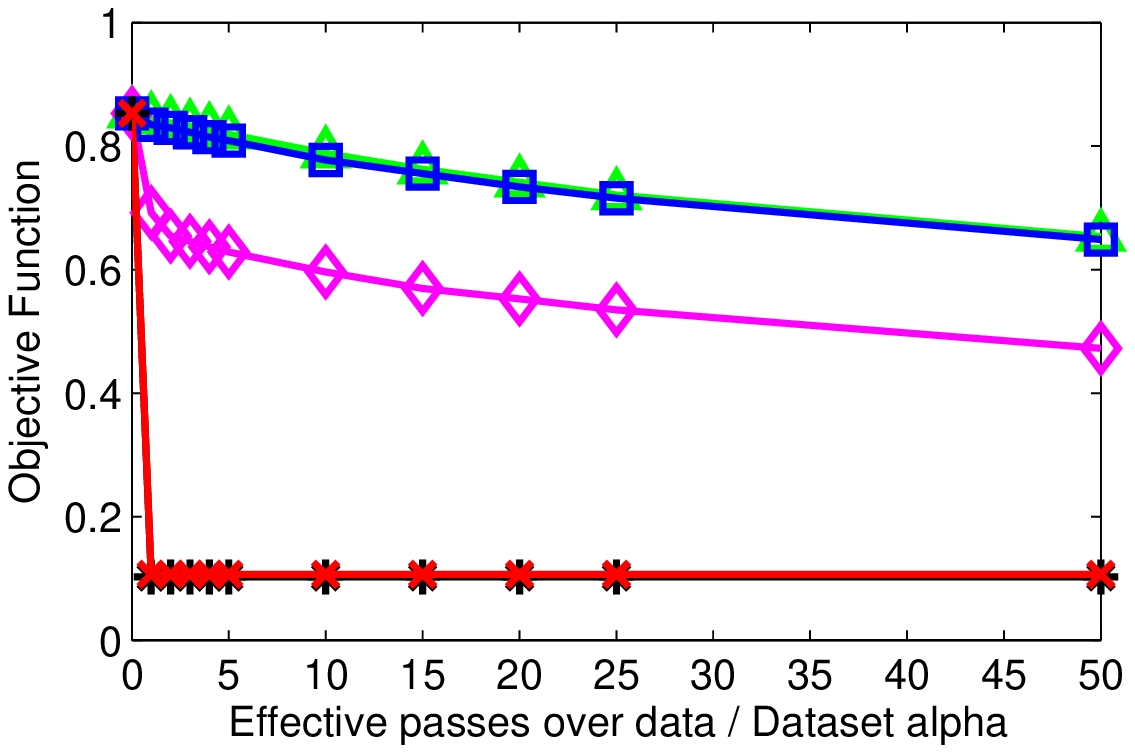}\hfill
   \includegraphics[width=0.33\linewidth]{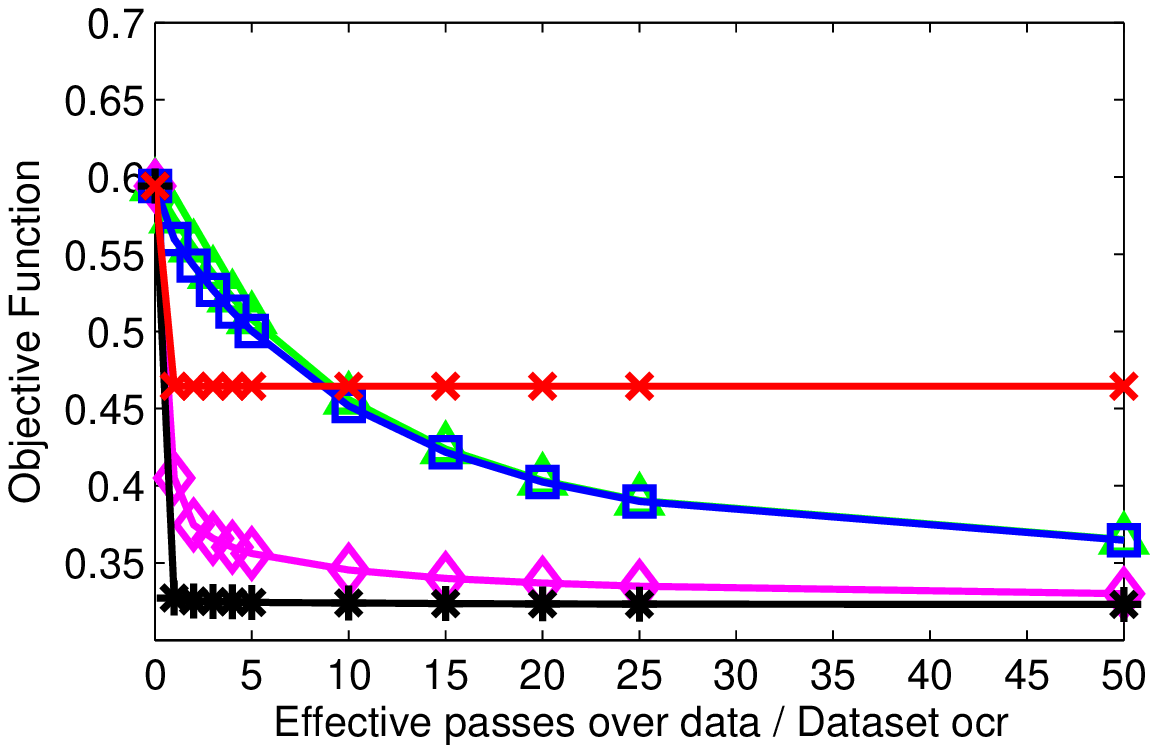}\\
   \includegraphics[width=0.33\linewidth]{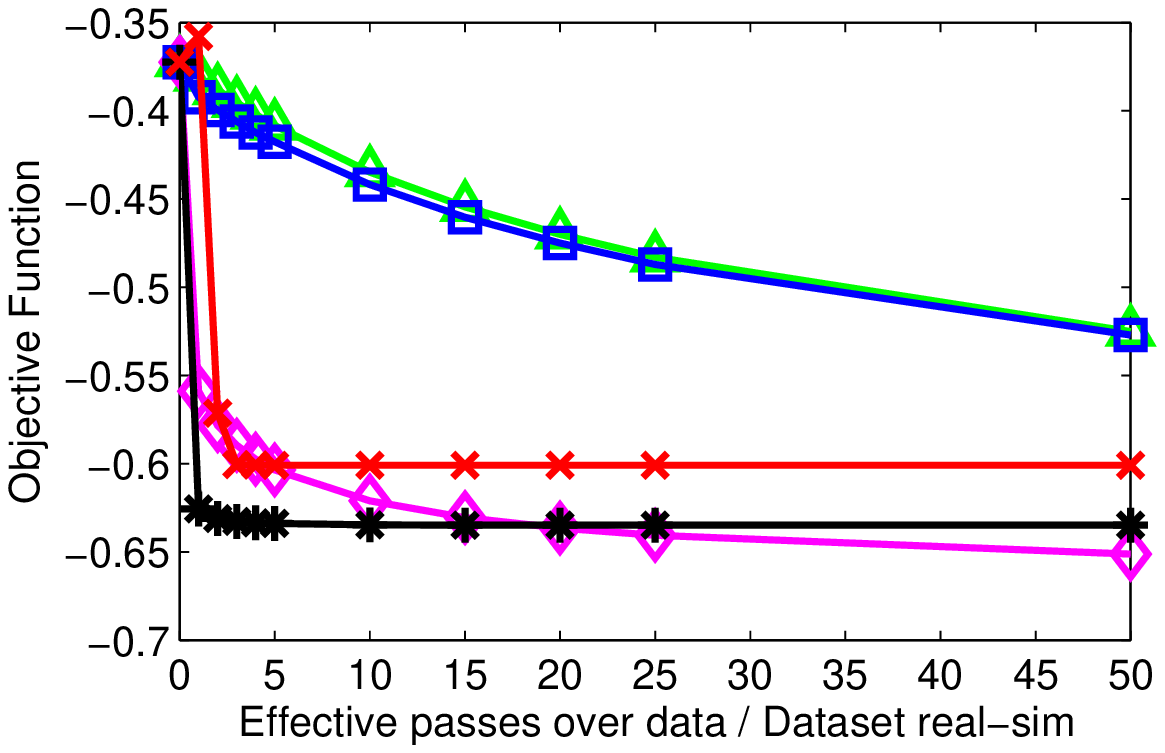}\hfill
   \includegraphics[width=0.33\linewidth]{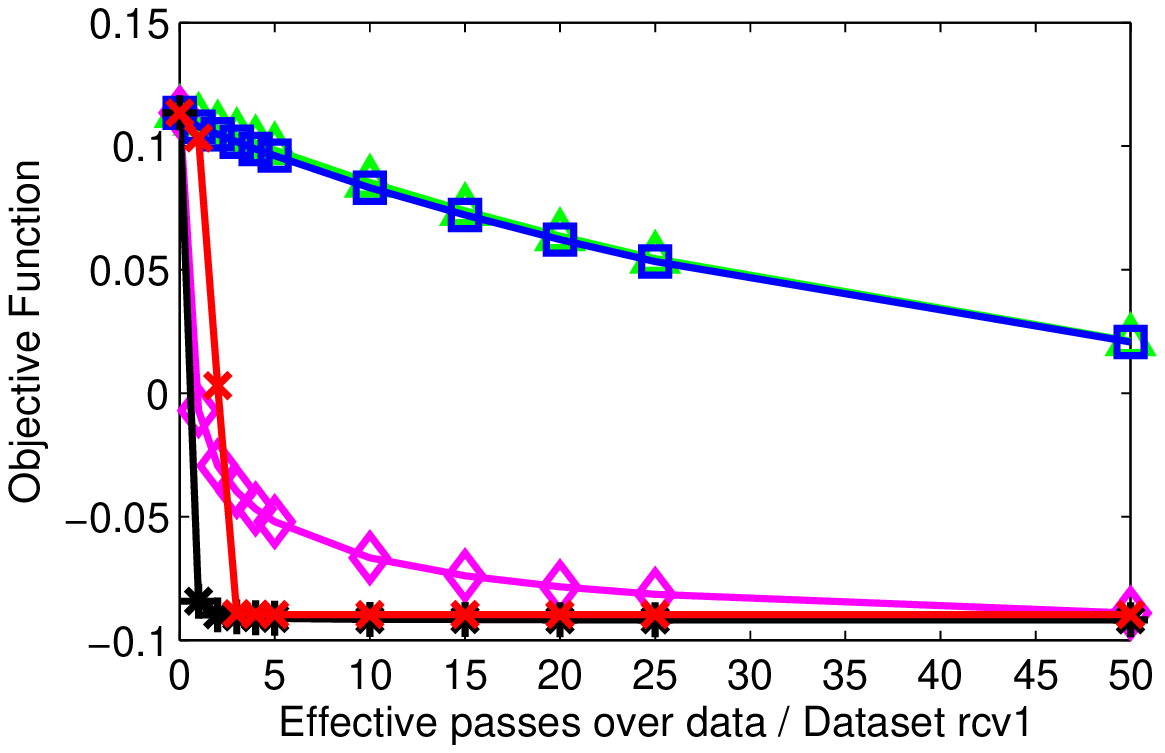}\hfill
   \includegraphics[width=0.33\linewidth]{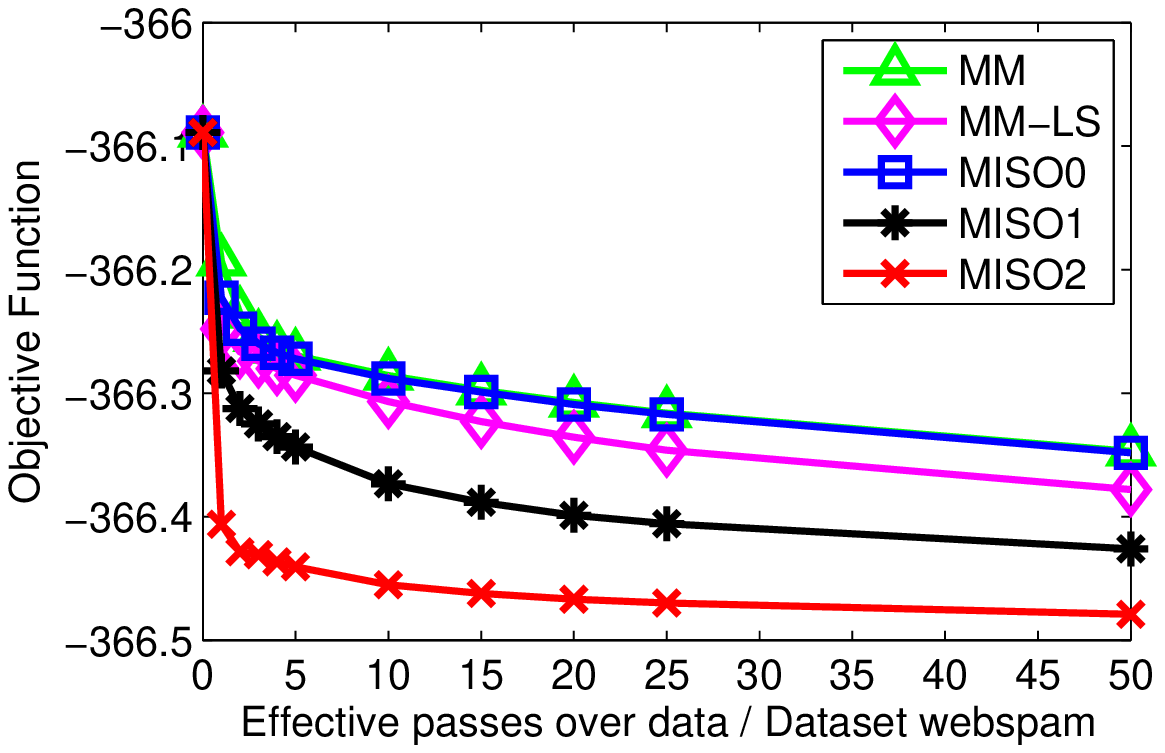}
   \caption{Objective function during the sparse estimation experiment.}
   \label{fig:dcepochs}
\end{figure}

\begin{figure}[hbtp]
   \centering
   \includegraphics[width=0.33\linewidth]{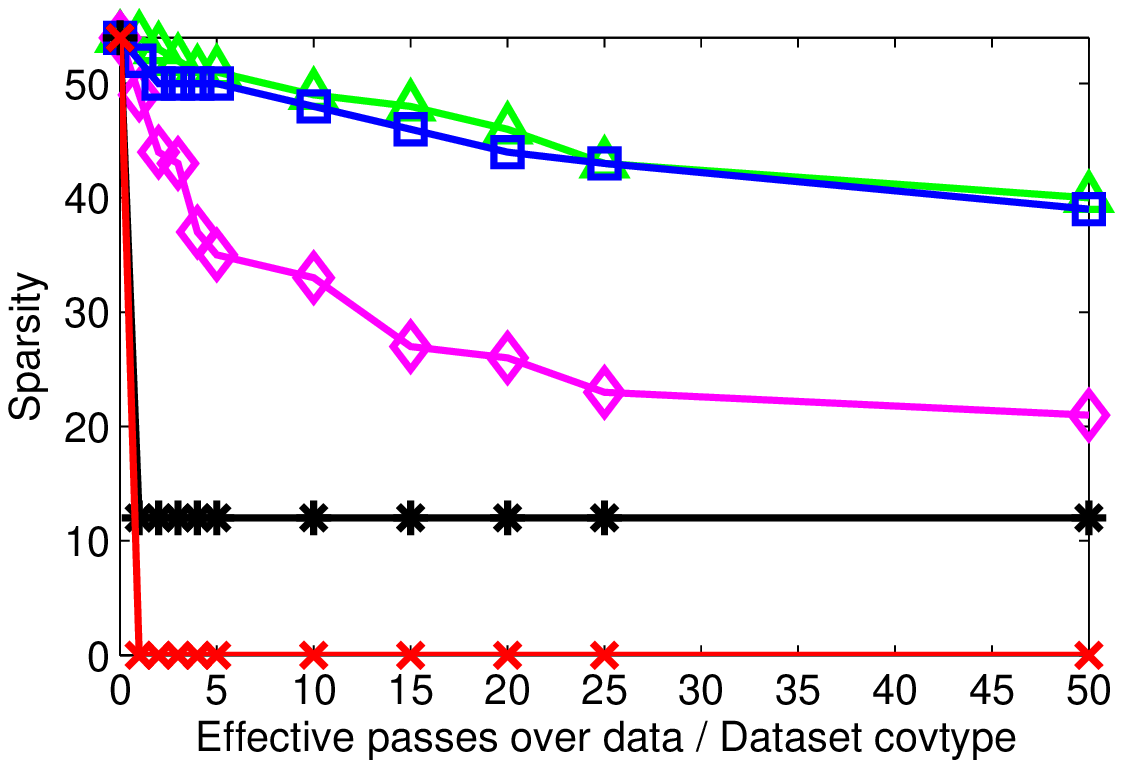}\hfill
   \includegraphics[width=0.33\linewidth]{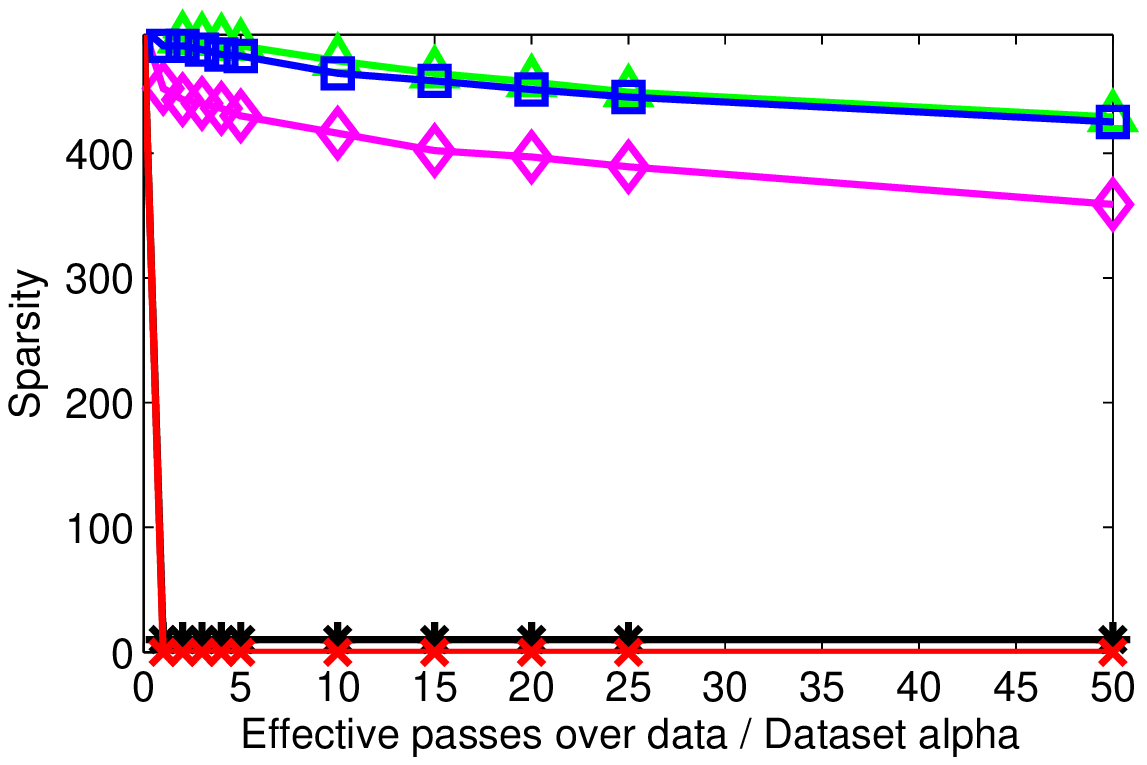}\hfill
   \includegraphics[width=0.33\linewidth]{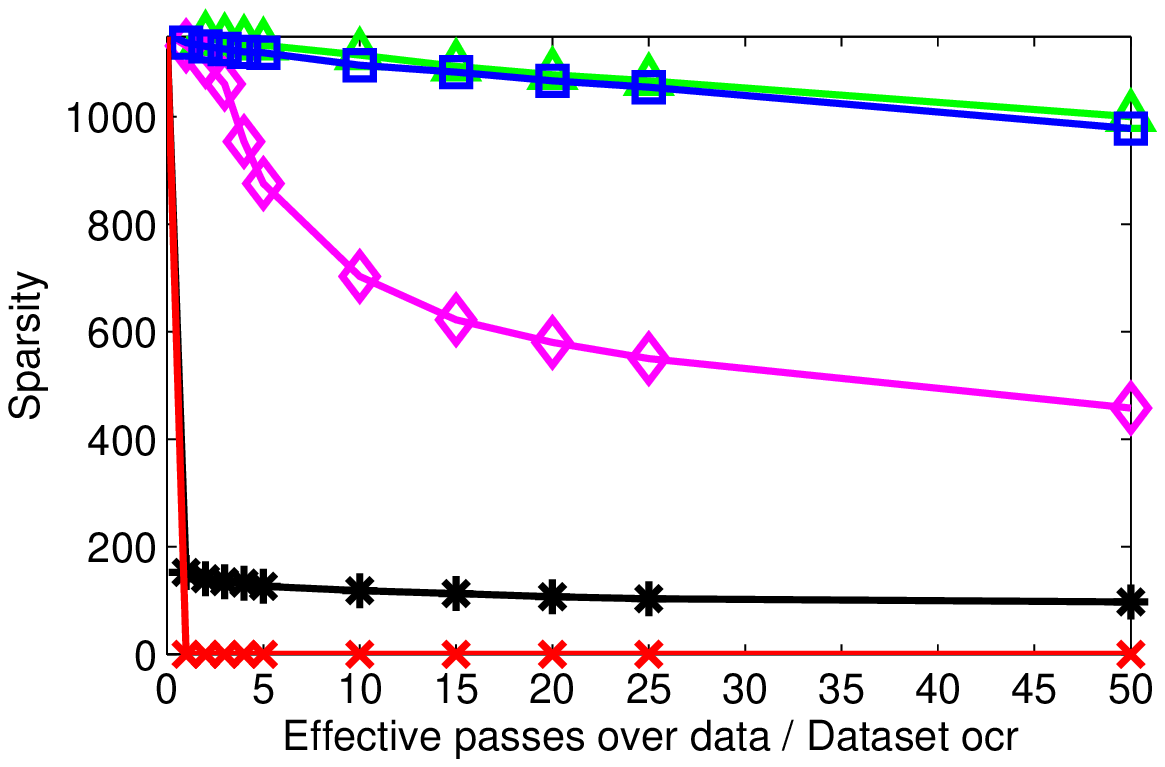}\\
   \includegraphics[width=0.33\linewidth]{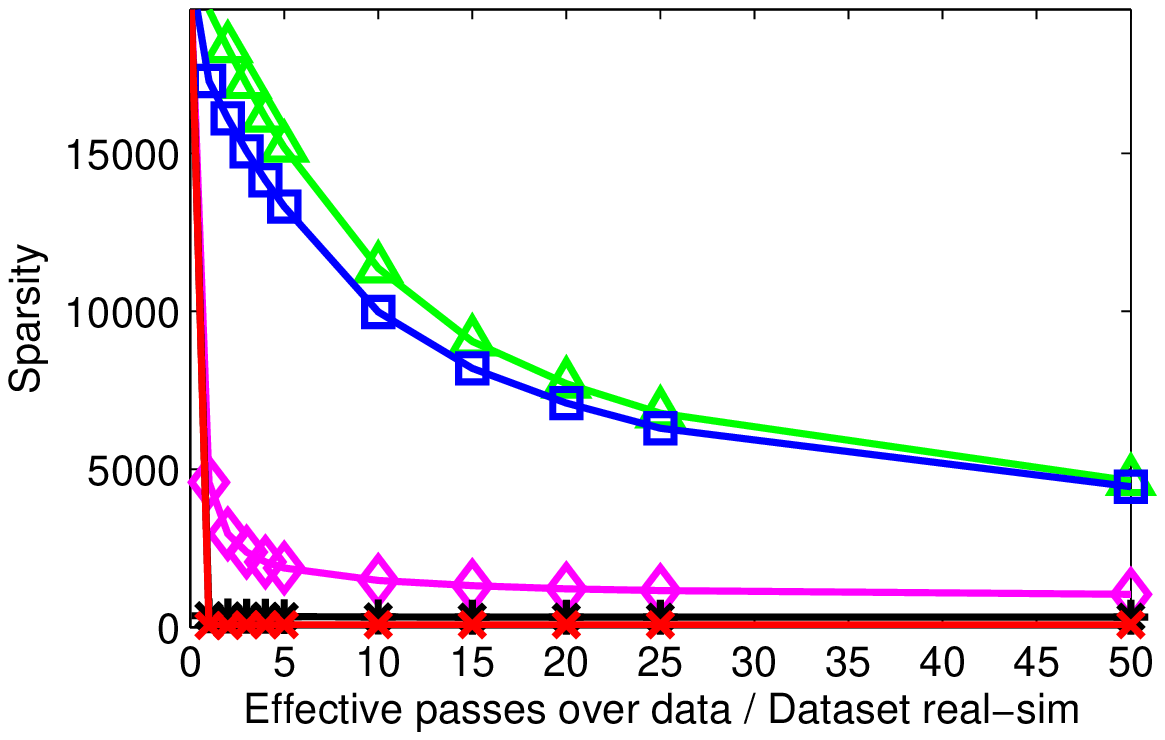}\hfill
   \includegraphics[width=0.33\linewidth]{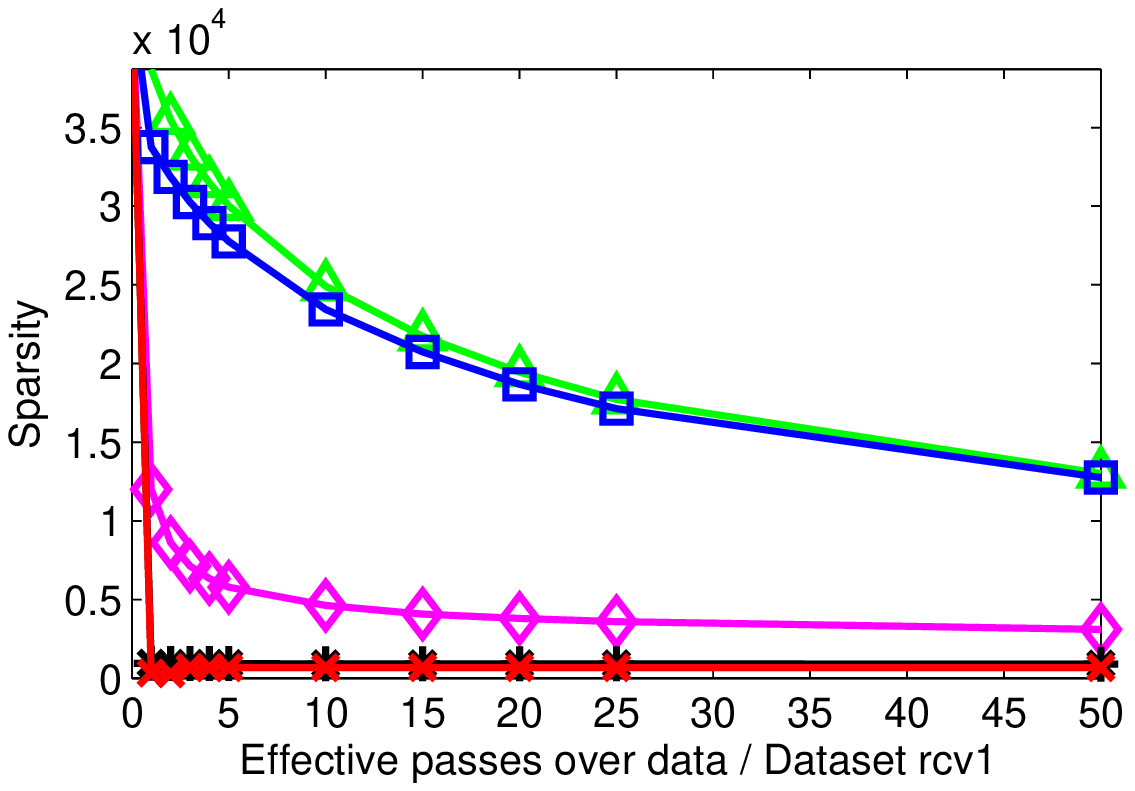}\hfill
   \includegraphics[width=0.33\linewidth]{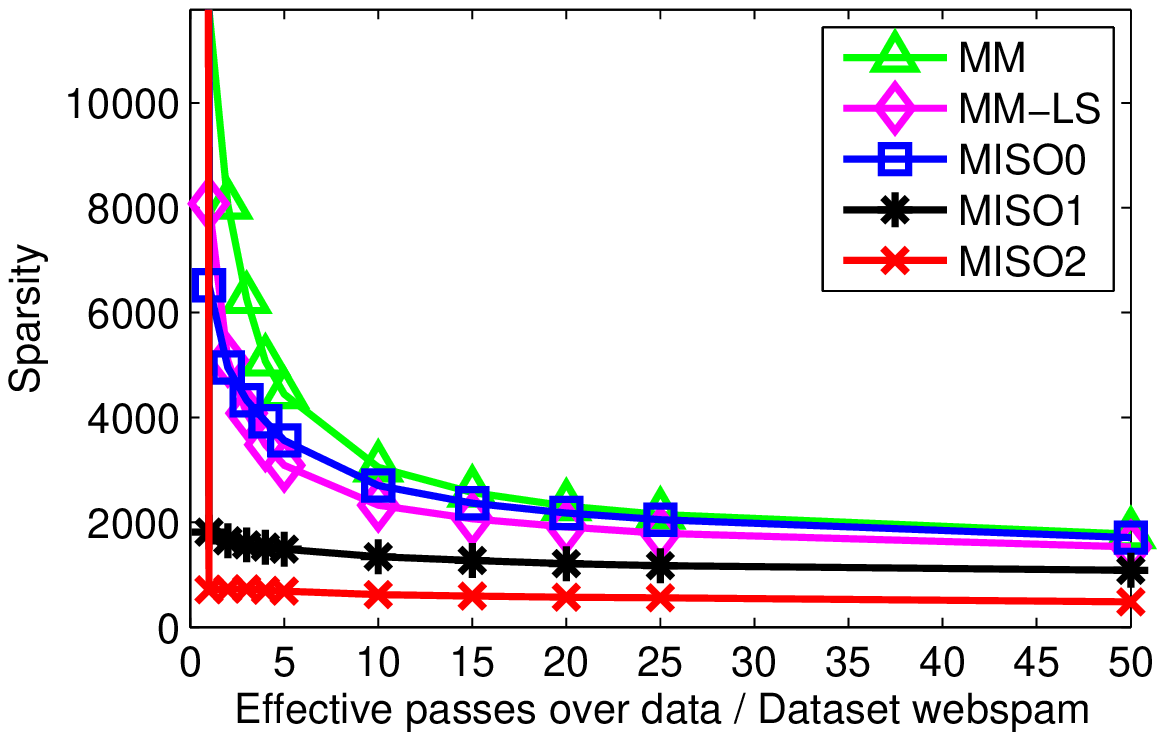}
   \caption{Sparsity of the solution during the sparse estimation experiment.}\label{fig:dcspars}
\end{figure}

The majorization-minimization principle is appealing for non-convex and
non-smooth optimization, where only few algorithms apply. Here, we address a
sparse estimation problem presented in Section~\ref{subsec:dc}:
\begin{equation}
   \min_{\theta \in \Real^p} \frac{1}{T} \sum_{t=1}^T \frac{1}{2}(y_t -
   \x_t^{\top}\theta)^2 + \lambda\sum_{j=1}^p \log(|\theta[j]|+\varepsilon),
   \label{eq:expdc}
\end{equation}
where the scalars $y_t$ and the vectors~$\x_t$ are the same as in the previous
section, and~$\varepsilon$ is set to $0.01$. The model parameter~$\lambda$
controls the sparsity of the solution. Even though~(\ref{eq:expdc}) is
non-convex and non-smooth, stationary points can be obtained in various ways.
In this section, we consider majorization-minimization approaches where the
penalty function $\theta \mapsto \sum_{j=1}^p \log(|\theta[j]|+\varepsilon)$ is
upper-bounded as in Eq.~(\ref{eq:upperbounddc}), whereas the functions $\theta
\mapsto (1/2)(y_t - \x_t^{\top}\theta)^2$ are upper-bounded by the Lipschitz
gradient surrogates of Section~\ref{subsubsec:gradient}. We compare five approaches:
\begin{itemize}
   \item[MM] Algorithm~\ref{alg:generic_batch} with the trivial Lipschitz constant $L=(1/T)\sum_{t=1}^T0.25\|\x_t\|_2^2$.
   \item[MM-LS] Algorithm~\ref{alg:generic_batch} with the line-search scheme of ISTA~\cite{beck} for adjusting~$L$.
   \item[MISO] we compare MISO0, MISO1, and MISO2, as in the previous section.
\end{itemize}
We choose a parameter $\lambda$ for each dataset, such that the solution with
the lowest objective function obtained by any of the tested method has
approximately a sparsity of~$10$ for datasets \textsf{covtype} and
\textsf{alpha}, 100 for \textsf{ocr} and \textsf{real-sim}, and $1\,000$ for
\textsf{rcv1} and \textsf{webspam}. The methods are initialized with $\theta_0= (\|\y\|_2/\|\X\X^\top \y\|_2)\X^\top \y$;
indeed, the initialization $\theta_0=0$ that was a natural choice in Section~\ref{subsec:explog} appears to be often a bad stationary point of problem~(\ref{eq:expdc}) and
thus an inappropriate initial point.
We report the objective function values
for different passes over the data in Figure~\ref{fig:dcepochs}, and the
sparsity of the solution in Figure~\ref{fig:dcspars}.
Our conclusions are the following:
\begin{itemize}
   \item methods with line searches do significantly better than those without,
      showing that adjusting the constant $L$ is important for these datasets;
   \item MISO1 does asymptotically better than MM-LS for five of the datasets
      after $50$ epochs and slightly worse for \textsf{real-sim}; in general,
      MISO1 seems to converge substantially faster than other approaches, both
      in terms of objective function and in terms of the support of the
      solution. 
   \item the performance of MISO2 is mitigated. In one
      case, it does better than MISO1, but in some others, it converges to the
      stationary point~$\theta\!=\!0$.
\end{itemize}

\section{Conclusion}\label{sec:ccl}
In this paper, we have presented new algorithms based on the
majorization-minimization principle for minimizing a large sum of functions.
The main asset of our approach is probably its applicability to a large class
of non-convex problems, including non-smooth ones, where we obtain convergence
and asymptotic stationary point guarantees.  For convex problems, we also
propose new incremental rules for composite optimization, which are
competitive with state-of-the-art solvers in the context of large-scale machine
learning problems such as logistic~regression.

We note that other majorization-minimization algorithms have recently been
analyzed, such as block coordinate variants in~\cite{mairal17,razaviyayn2} and
stochastic ones in~\cite{choromanska,mairal18,razaviyayn}. In particular, we have
proposed in~\cite{mairal18} a stochastic majorization-minimization algorithm that does not require to
store information about past iterates, when the objective function is an
expectation. Since the first version of our work was published
in~\cite{mairal18}, MISO has also been extended by other authors
in~\cite{zhong} using the alternating direction method of multipliers
framework.

For future work, we are currently investigating extensions of the scheme
MISO$\mu$ for strongly convex objective functions. We believe that the
algorithm can be modified to remove the large sample condition~$T \geq 2L/\mu$,
that the convergence proof can be extended to the proximal setting, and that it
is possible to use acceleration techniques in the sense of
Nesterov~\citep{nesterov4}.  Another interesting direction of research would be
to study the stability of our result to inexact minimization of surrogate
functions following for instance the analysis of~\cite{schmidt} for proximal
gradient methods.

\section*{Acknowledgments}
The author would like to thank Zaid Harchaoui, Francis Bach, Simon
Lacoste-Julien, Mark Schmidt, Martin Jaggi, the associate editor, and the
anonymous reviewers for their useful comments.

\appendix
\section{Basic definitions and useful results}\label{appendix:background}
The following definitions can be found in classical 
textbooks, e.g,~\cite{bertsekas,borwein,boyd,nocedal}.  For the sake of
completeness, we briefly introduce them here.
\begin{definition}[{\rm Directional derivative}]\label{def:derivative}
   Let us consider a function $f: \Real^p \to \Real$ and $\theta, \theta'$ be in~$\Real^p$. When it exists,
   the following limit is called the directional derivative of $f$ at $\theta$ in
   the direction $\theta'-\theta$:
   $ \nabla f(\theta,\theta'-\theta) \defin \lim_{t \to 0^+}
   {(f(\theta+t(\theta'-\theta)) - f(\theta)})/{t}.$
   When $f$ is differentiable at $\theta$, 
   directional derivatives exist in every direction, and $\nabla
   f(\theta,\theta'-\theta)=\nabla f(\theta)^\top (\theta'-\theta)$.
\end{definition}

\begin{definition}[{\rm Stationary point}]\label{def:stationary}
   Let us consider a function $f: \Theta \subseteq \Real^p \to \Real$, where $\Theta$ is a convex set, such that $f$ admits a directional derivative $\nabla f(\theta,\theta'-\theta)$ for all $\theta,\theta'$ in $\Theta$. We say that $\theta$ in~$\Theta$ is a stationary point if for all $\theta'$ in~$\Theta$,
   $  \nabla f(\theta, \theta'-\theta) \geq 0$.
\end{definition}
\begin{definition}[{\rm Lipschitz continuity}]\label{def:lipschitz}
   A function $f: \Real^p \to \Real$ is $L$-Lipschitz continuous for some $L > 0$ when for all $\theta,\theta'$ in~$\Real^p$,
   $      |f(\theta')-f(\theta)| \leq L \|\theta-\theta'\|_2.$
\end{definition}
\begin{definition}[{\rm Strong convexity}]\label{def:strong_convexity}
   Let $\Theta$ be a convex set. A function $f: \Theta \subseteq \Real^p \to \Real$ is called $\mu$-strongly convex when there exists a constant $\mu > 0$ such that for all $\theta'$ in $\Theta$, the function $\theta \mapsto f(\theta)-\frac{\mu}{2}\|\theta-\theta'\|_2^2$ is convex. 
\end{definition}

We now present two lemmas that are useful for characterizing first-order
surrogate functions. Their proofs can be found in the appendix
of~\cite{mairal17}.
\begin{lemma}[\rm Regularity of residual functions]\label{lemma:convexerror}
   Let $f,g: \Real^p \to \Real$ be two functions. Define $h\defin g-f$. Then,
      if $g$ is $\rho$-strongly convex and $f$ is $L$-smooth, with $\rho \!\geq\! L$, $h$ is $(\rho\!-\!L)$-strongly convex; \label{step:convexerror}
       if $g$ and $f$ are convex and $L$-smooth, $h$ is also $L$-smooth;
       if $g$ and $f$ are $\mu$-strongly convex and $L$-smooth, $h$ is $(L-\mu)$-smooth.
\end{lemma}
\begin{lemma}[\rm Regularity of optimal value functions]\label{lemma:danskin}
   Let $f: \Real^{p_1} \times \Theta_2 \to \Real$ be a function of two
   variables where $\Theta_2 \subseteq \Real^{p_2}$ is a convex set. Assume that
   \begin{itemize}
      \item $\theta_1 \mapsto f(\theta_1,\theta_2)$ is differentiable for all~$\theta_2$ in $\Theta_2$;
      \item $\theta_2 \mapsto \nabla_{1} f(\theta_1,\theta_2)$ is $L$-Lipschitz continuous for all $\theta_1$ in $\Real^{p_1}$;
      \item $\theta_2 \mapsto f(\theta_1,\theta_2)$ is $\mu$-strongly convex for all $\theta_1$ in~$\Real^{p_1}$. 
   \end{itemize}
   Also define $\tilde{f}(\theta_1) \defin
   \min_{\theta_2 \in \Theta_2} f(\theta_1,\theta_2)$.
   Then, $\tilde{f}$ is differentiable and $\nabla \tilde{f} (\theta_1) = \nabla_{1} f(\theta_1,\theta_2^\star)$, where
   $\theta_2^\star \defin \argmin_{\theta_2 \in \Theta_2} f(\theta_1,\theta_2)$.
 Moreover, if $\theta_1 \mapsto \nabla_{1} f(\theta_1,\theta_2)$ is $L'$-Lipschitz continuous for all $\theta_1$ in $\Real^{p_1}$,
 the gradient~$\nabla \tilde{f}$ is $(L'+L^2/\mu)$-Lipschitz.
\end{lemma}

\bibliographystyle{siam}

\end{document}